\documentclass[12pt]{article}

\usepackage{amssymb}
\usepackage{latexsym}

\newtheorem{theorem}{Theorem}
\newtheorem{lemma}[theorem]{Lemma}

\newtheorem{proposition}[theorem]{Proposition}
\newtheorem{corollary}[theorem]{Corollary}

\newtheorem{conjecture}{Conjecture}

\newenvironment{proof}{\paragraph{\it Proof.}}{$\square$ \vskip0.4cm}
\newenvironment{remark}{\paragraph{\it Remark.}}{\vskip0.4cm}

\newcommand{\Z}{{\mathbb Z}}
\newcommand{\R}{{\mathbb R}}
\newcommand{\Q}{{\mathbb H}}

\begin{document}


\title{Inscribing cubes and covering by rhombic dodecahedra via
equivariant topology}

\author{T. Hausel \thanks{\noindent 
Research (partially) supported by Hungarian National Foundation for
Scientific Research, grant no. 23444, and by Trinity College,
Cambridge, UK}\and E. Makai, Jr.
\thanks{Research (partially) supported by Hungarian National Foundation for 
Scientific Research, grant no. T-030012}\and  A. Sz\H ucs\thanks{Research 
(partially) supported by Hungarian National Foundation for 
Scientific Research, grant no. A 046/96}}

\maketitle

\begin{abstract} 
First, we prove a special case of Knaster's problem, implying
that each symmetric convex body in $\R^3$ admits an inscribed 
cube. We deduce it from a theorem in equivariant topology, which says
that there is no $S_4$-equivariant map from $SO(3)$ to $S^2$, where
$S_4$ acts on $SO(3)$ as the rotation group of the cube and on $S^2$
as the symmetry group of the regular tetrahedron. 
We also give some generalizations.

Second, we show how the above non-existence theorem yields Makeev's
 conjecture
in $\R^3$ that each set
in $\R^3$ of diameter $1$ can be covered by a rhombic dodecahedron, 
which has distance $1$ between its opposite faces. This reveals an
unexpected connection between inscribing cubes into symmetric bodies
and covering sets by rhombic dodecahedra. 
  
Finally, we point out a possible application of our second theorem to 
the Borsuk problem in $\R ^3$.

\vskip.4cm
\noindent
{\bf Keywords and phrases:} convex bodies, inscribed polyhedra,
universal covers,
cubes, boxes, Knaster's problem, 
equivariant algebraic topology, 
sphere bundles, vector bundles, intersection numbers, obstruction class, 
Stiefel-Whitney class,  
\vskip.4cm
\noindent
{\bf 1991 Mathematics Subject Classification.} Primary: 52A15; Secondary:
55Mxx
\end{abstract}
\newpage


\section{Introduction}
\label{introduction}

Problems of inscribing and circumscribing 
polyhedra to convex bodies have a vast
literature; for a comprehensive recent survey cf. in \cite{klee-wagon}, Ch. 11;
 cf. also \cite{makeev-1}.
In \cite{hausel-makai-szucs} we have surveyed some of the
literature and announced some of the results of the present paper.

The proofs of 
such theorems are frequently based on theorems for continuous maps,
such as some special cases of the so called Knaster problem.  

The Knaster problem asks,  
for a given continuous function $F:S^{n-1}\rightarrow \R^k$ 
and a finite 
subset $X=\{x_1,...,x_m\}$ of $S^{n-1}$, whether there is a rotation 
$A\in SO(n)$ such that $F(Ax_1)=\dots=F(Ax_m)$. 
Observe that, choosing $F$ linear, for a positive answer it
is necessary that $X$ should lie in an $(n-k)$-plane; correspondingly, this 
question is frequently asked for $m=n-k+1$ only, as originally in
\cite{knaster}. The answer to the Knaster problem in general is
negative: counterexamples were given by Makeev in \cite{makeev1} and
\cite{makeev2}, by Babenko and Bogatyi in \cite{babenko-bogatyi} 
and recently by Chen in \cite{chen}.

However there are many positive results: a famous special case
is the Borsuk-Ulam theorem which states that for every continuous map
$F:S^{n-1}\rightarrow \R^{n-1}$  there is an $x\in S^{n-1}$
such that $F(x)=F(-x)$. This is Knaster's problem when $k=n-1$, $m=2$
and $X=\{e_1,-e_1\}$. It was generalized by Hopf in \cite{hopf} for any
two points in $S^{n-1}$. Other positive results are proved when $k=1$  and 
$X$ is the set of vertices of any regular $(n-1)$-simplex
inscribed to $S^{n-1}$ in \cite{bogatyi-khimshiashvili},
$X=\{e_1,\dots,e_{n}\}$ \cite{yamabe-yujobo},
$X=\{ \pm e_1,...,\pm e_{n-1}\}$ Corollary 7.5, \cite{yang},
 $X$ any $3$-element subset 
of $S^2$ \cite{floyd}, $X=\{ \pm x, \pm y \} \subset S^2$ 
\cite{dyson}, \cite{livesay}, $X$ being 
the vertex set of any planar rectangle inscribed to $S^2$
Theorem D, \cite{griffiths}\footnote{ We note, that a slight criticism of the 
Griffiths paper appeared in 
\cite{makeev5}. This states that another theorem of [Griff],
namely his Theorem C is correctly proved by calculating intersection
numbers mod 2, not in $\Z$ as in [Griff]. Since [Griff] uses large
computer calculations, we are unable to decide in this point. [Griff],
Theorem C is not used in proving [Griff], Theorem D.}. 

In this paper we settle the following

\begin{theorem}
\label{cube} 
Let $F:S^2 \to {\R}$ be a continuous even function, 
i.e., $F(x)=F(-x)$, and let $\{\pm
v_i \mid 1 \le i \le 4 \}$ be the vertices of a cube inscribed to $S^2$. Then
there exists an $A \in SO(3)$ such that $F(\pm Av_1)=F(\pm Av_2)=F(\pm Av_3)=
F(\pm Av_4)$.
\end{theorem}

Of course, this theorem also follows from Theorem D of \cite{griffiths}
cited 
above. However the
proof of that theorem is quite involved, and we hope our proof of this special 
case is considerably simpler. Our proof avoids calculations, and uses only 
elements of equivariant algebraic topology. Hopefully our theorems
about non-existence of certain equivariant maps are of independent
interest.  

Note that when $F$ is the Minkowski norm associated to a symmetric convex body
in $\R^3$,
Theorem~\ref{cube} proves that each symmetric convex body in $\R^3$ contains an
inscribed cube. 

In Section~\ref{konkavfejezet} 
we will prove a more involved 
generalization of Theorem~\ref{cube} which is no longer a 
consequence of Theorem D of \cite{griffiths}. It says:

\begin{theorem}
\label{konkav} 
Let $f:S^2 \to \R ^3$ be a $C^1$ embedding of $S^2$ in $\R ^3$, satisfying 
$f(-u)=-f(u)$ for $u \in S^2$. Further suppose that $f$ is homotopic 
to the standard embedding, 
via a $C^1$ homotopy $H:S^2 \times [0,1] \to \R^3 \setminus \{0\}$ satisfying 
$H(-u,t)=-H(u,t)$ for $u \in S^2$ and $t \in [0,1]$.
Then there is a cube in 
$\R^3$, with its centre in
the origin, having all its vertices on the surface $f(S^2).$
\end{theorem}

The second theme of the present paper concerns {\em universal covers}.
One says that a set $C \subset {\R}^n$ is a universal cover 
in ${\R}^n$ 
(cf. \S 10, 47, p. 87, \cite{bonnesen-fenchel} and \cite{makeev3}),
if for any set $X \subset {\R}^n$ of diameter at most $1$ there exist
$A \in O(n)$ and $x \in {\R}^n$, such that $X \subset AC+x$. 
(If $C$ has a plane of symmetry, we can write here $A \in SO(n)$. In any case,
since $O(n)$ is disconnected, when proving the universal cover property of 
some $C \subset \R^n$, it is
 reasonable to use only $A \in SO(n)$.)  
Since the 
diameter of $X$ and that of its closed convex hull are equal, we may assume 
that $X$ is compact, convex, non-empty. For example, the unit ball is a 
universal cover, but also the circumball of a regular simplex with unit 
edges is a universal cover [Jung 1901], [Jung 1910], cf. also 
\cite{bonnesen-fenchel}, \S 10, 44, (7)
(this ball has radius $\sqrt{n/(2n+2)}$).

Moreover P\'al in \cite{pal} showed that a regular hexagon, with distance of 
opposite sides equal to 
$1$, was also a universal cover in $\R^2$. It is an interesting open 
question, due to Makeev (\cite{makeev3}, p.129; cf. also our 
Section~\ref{covers}),
whether an analogue of this holds in ${\R}^n$.

Here we prove Makeev's conjecture in $\R^3$:

\begin{theorem} A rhombic dodecahedron $U_3$, with distance of opposite
faces equal to $1$, is a universal cover in $\R^3$. I.e., any set of
diameter $1$ in $\R^3$ can be covered by $U_3$. 
\label{rhombic}
\end{theorem}

We prove both Theorem~\ref{cube} in Section~\ref{cubes}
and Theorem~\ref{rhombic} in Section~\ref{covers},
as a consequence of the following result in equivariant topology:

\begin{proposition}
\label{s4}
There is no equivariant map 
$$g:(SO(3),\rho_{S_4})\rightarrow (S^2,\tau_{S_4}),$$
where $\rho_{S_4}$
is the $S_4$ action on $SO(3)$ 
given by multiplication from right by the group of 
rotations of the cube: $S_4\subset SO(3)$, and $\tau_{S_4}$ is the
$S_4$ action on $S^2$ given by the symmetry group of a regular
tetrahedron inscribed in $S^2$.   

\end{proposition}

Following \cite{jerrard} and \cite{griffiths},
the proof is based on the idea of {\em test functions}. 
Namely we show
that one can inscribe exactly one cube in an ellipsoid $E$ of different
lengths of axes. From this in turn we get a test function 
$f_E:(SO(3),\rho_{S_4})\rightarrow (\R^3,\tau_{S_4})$ which has
exactly one transversal zero, which will prove the proposition.  

In connection with 
this proposition, 
we investigate the existence of equivariant maps $SO(3) \to S^2$ for 
certain group actions of subgroups of $S_4$
on $SO(3)$ and $S^2$ in Section~\ref{equivariant}. 
As a consequence we
deduce Theorem~\ref{square} in Section~\ref{squares},
which is Theorem~\ref{cube} with the cube replaced by a square
based box (a box is a rectangular parallelepiped).
Thus one can inscribe a square based box, similar to any
given one,
into any symmetric
convex body in $\R^3$.

In Section~\ref{box} we show that Theorem D in \cite{griffiths} easily
yields a generalization of Theorem~\ref{cube}, 
namely the same statement with
the vertices of the cube replaced by those of any box. Consequently 
each symmetric convex body in $\R^3$ admits an inscribed box similar to a given
one. 

In Section~\ref{konkavfejezet}
we prove Theorems~\ref{konkav} and ~\ref{konkav'}, that generalize
Theorems~\ref{cube} and ~\ref{square}, showing
that certain centrally symmetric surfaces in $\R ^3$ admit inscribed
cubes, and inscribed similar copies of any square based box,
respectively.

In the final Section~\ref{covers} we prove Theorem~\ref{rhombic} as a 
consequence  of Proposition~\ref{s4}. We finish the paper 
by pointing out a possible application of Theorem~\ref{rhombic} to
Borsuk's problem in $\R^3$, which states that any set of diameter $1$
in $\R^3$ can be divided into $4$ sets of diameter smaller than $1$.

We have announced the main results of this paper in
\cite{hausel-makai-szucs} 
in 1997. Since 
the present paper was written 
we learnt of the paper \cite{makeev4} by 
V. V. Makeev from 1997 and of 
\cite{kuperberg} by G. Kuperberg from 1998. Both papers consider similar
questions. In \cite{kuperberg} 
Proposition~\ref{s4} is also proved and applied to problems similar 
to ours.


\vskip.4cm 

\paragraph{\bf Acknowledgement.} 
The authors express their gratitude to M. Lassak for 
having turned their attention to the question of inscribing cubes
into symmetric bodies
(\cite{lassak}), when the
second named author enjoyed his hospitality in Bydgoszcz in 1992, and
to V. Klee for useful e-mail correspondence.

\newpage


\section{Inscribed cubes: proof of Theorem~\ref{cube}}
\label{cubes}

Eggleston in \cite{eggleston} has
asserted that there exists a centrally symmetrical convex body
in ${\R}^3$ that does not admit an inscribed cube. For this he has stated a 
lemma (p.79). Formula (2) in the proof 
of this lemma should stand correctly 
$$x^2/A^2+y^2/B^2+z^2/C^2+\lambda (xy+yz+zx)=d^2/A^2+d^2/B^2+d^2/C^2
-\lambda d^2.$$ In fact, the statement of the lemma itself is not 
correct. Namely, it asserts the following. Let 
$$V=\{ (x,y,z) \in 
\R ^3 \mid x,y,z \in \{ -1,1\},\,\, x+y+z \in \{-1,1\}\}.$$ Suppose 
that the boundary of an ellipsoid $E \subset \R ^3$ contains $V$. Then the 
lemma asserts that either $E$ has two equal axes, or else its 
centre is $0$, its axes 
are parallel to the coordinate axes, and its boundary also contains 
$(-1,-1,-1),\,\,(1,1,1)$.

However, $V$ consists of the 
vertices of an octahedron (i.e., an affine image of a regular 
octahedron, 
i.e., $V= \{ \pm x, \pm y, \pm z\}$, with $x,y,z \in \R^3$ 
linearly independent). Further, any ellipsoid has lots of inscribed 
octahedra, with no connection to the directions of the axes of the 
ellipsoid. 

To have a concrete counterexample, let us consider an 
ellipsoid $E_0$ of centre $0$, of equation $0.5x^2+y^2+1.5z^2=3$, 
whose boundary 
passes through the points $(\pm 1,\pm 1,\pm 1)$.
Its section with the plane $z=1$ has equation $z=1,\,\,x^2/3+2y^2/3=1$. 
Then let $\epsilon>0$ be 
sufficiently small, and let us consider a conic $C$ in the plane $z=1$ 
passing through the five points $(-1,-1,1),\,\,(-1,1,1),\,\,(1,-1,1), 
\,\, (1+ \epsilon, 1, 1), (\sqrt 3, 0,1)$. This is uniquely 
determined, depends continuously on $\epsilon$, is an ellipse, and 
does not pass through $(1,1,1)$, since 
it intersects the line $y=z=1$ in two points only. Let the equation 
of $C$ be $z=1,\,\,\sum _{i,j} a_{ij} x_i x_j +\sum _i b_i x_i +1 
=0\,\,(a_{ij}=a_{ji})$, where $x_1=x,\,\,x_2=y$, and the 
coefficients depend continuously on $\epsilon$. Let $E$ be an 
ellipsoid of centre $0$, 
with $\partial E \supset C$, of equation $\sum _{i,j} A_{ij} x_i x_j +C 
=0\,\,(A_{ij}=A_{ji})$, with $x_3=z$. We may suppose $a_{ij} 
=A_{ij} \,\,(i,j \le 2), \,\, b_i=2A_{i3}\,\,(i \le 
2),\,\,1=A_{33} +C$.  Let now $A_{33}=1.5$; this determines 
the equation of $E$ uniquely, and the coefficients in this equation 
depend on $\epsilon$ continuously. For $\epsilon =0$ we have the 
original ellipsoid $E_0$. 

Therefore for $\epsilon$ sufficiently small, 
$E$ is an ellipsoid, and the corresponding quadratic form  
has three different eigenvalues, hence $E$ has three different axis lengths. 
By construction, $\partial E \supset V, \,\, (1,1,1) \not\in \partial E$, 
contradicting the claim of the lemma of \cite{eggleston}.

\vskip.4cm

We will actually prove the opposite of Eggleston's statement as a
corollary to Theorem~\ref{cube}: 

\begin{corollary}
Any centrally symmetric convex body $K$ in $\R^3$ possesses
an inscribed cube, centred at the centre of $K$.
\label{inscribe}
\end{corollary}

Before we proceed to the proof 
 we recall the symmetry group of the regular tetrahedron
and the rotation group (i.e., the group of orientation preserving 
symmetries) of the cube. 

The symmetry group of the regular tetrahedron (which we think of as 
inscribed to $S^2$) 
is $S_4$, the symmetric group on four letters. 
It is because the symmetries of a regular tetrahedron are in bijective
correspondence with the permutations of its vertices. Moreover 
the rotation group 
of the regular tetrahedron is clearly $A_4$, the alternating group on
four letters. 

The rotation group of the cube (which we also think of as inscribed to 
$S^2$)
is $S_4$. It can be seen by considering the
diagonals of the cube as follows. The cube has four diagonals and every 
rotation gives a permutation of these. Only the identity rotation maps
each diagonal to itself. 
On the other hand one can choose a vertex of each diagonal
yielding a regular tetrahedron of edge length $\sqrt 2$. Every permutation of the diagonals 
gives rise to a permutation of the vertices of the tetrahedron, which 
in turn determines a symmetry of the tetrahedron and thus that of the cube.
If it was not a rotation one could compose it 
with the antipodal map, yielding
a rotation of the cube inducing the given permutation of the diagonals. 
We will think of this rotation group as $S_4\subset SO(3)$ a subgroup
of $SO(3)$.    

Now, for the reader's convenience, we describe this embedding of $S_4$ 
to $SO(3)$. If $g \in S_4$ acts on the vertices of the regular tetrahedron as
an even permutation (i.e., in disjoint cycle representation has one of 
the forms: 
identity, $(ij)(kl)$ or $(ijk)$), then its embedded image $\iota (g)$ in 
$SO(3)$ is that element of $SO(3)$, that is the extension of 
the permutation $g 
\in S_4$ of the vertices of the tetrahedron. Geometrically, if the cube has 
vertices $\{ (\pm 1/\sqrt 3, \pm 1/\sqrt 3, \pm 1/\sqrt 3)\}$, 
these rotations are the identity, the 
rotations through $\pi$ about the coordinate axes, and the rotations through
$\pm 2\pi /3$ about the spatial diagonals of the cube, respectively. If $g \in
S_4$ acts on the vertices of the regular tetrahedron as an odd permutation 
(i.e., has one of the forms
$(ij)$ or $(ijkl)$), then $\iota (g)$ is $-1$ times
that element of $O(3)$, that is the extension of the permutation $g \in S_4$.
Geometrically, these rotations are the rotations through $\pi$ about all angle
bisectors of all pairs of the coordinate axes, and the rotations through $\pm
\pi /2$ about the coordinate axes, respectively. From now on, we will
write $g$ for
$\iota (g)$.   

\paragraph{\it Proof of Corollary~\ref{inscribe}.}
Let $C\subset \R^3$ be a fixed cube centred at the origin 
with diagonals $(v_i,-v_i)$, where $\|v_i\|=1$. We suppose that the centre of
$K$ is the origin, and 
let $\|\,\,\|_K$ be the Minkowski norm in $\R^3$, associated to $K$. 
Define 
the map $f: SO(3)\rightarrow \R^4$ sending $A\in SO(3)$ to 
$(\|A(v_1)\|_K,\|A(v_2)\|_K,\|A(v_3)\|_K,\|A(v_4)\|_K)$. 
Our task now is to show that the image
of $f$ intersects the diagonal 
$\Delta:=\{ (\lambda,\lambda,\lambda,\lambda):\lambda\in \R\}$ of
$\R^4$ since if $f(A)=(\lambda,\lambda,\lambda,\lambda)\,\,(\lambda >0)$ 
is a point on the
diagonal then the cube $\lambda^{-1}A(C)$ is inscribed in $K$, and is centred 
at the centre of $K$. 

The function $F=\|\,\,\|_K\mid_{S^2}:S^2\rightarrow\R$ is an even function
because $K$ is symmetric. Thus Corollary~\ref{inscribe} is a corollary
of Theorem~\ref{cube} indeed. $\square$

To prove Theorem~\ref{cube} we begin with two lemmas. 
The first one is surely well-known; we include its
short proof for the reader's convenience. Its statement and proof are slight 
variants of those of Theorem~1 of \cite{duncan-khavinson-shapiro}.
A special case of its corollary is contained in \cite{jerrard}, p. 240, proof 
of Lemma 6.

\begin{lemma}
Let $P \subset \R^n$ be a box, and let $E \subset \R^n$
be an ellipsoid circumscribed about $P$. Then $E$ has centre at the centre of 
$P$, and it has a system of $n$ axes parallel to the edges of $P$.
\label{center}
\end{lemma}

\begin{proof} Let the vertices of $P$ be $\{(\pm d_1,..., \pm d_n)\}$. Let 
$\partial E = \{(x_1,...,x_n) \mid  F(x_1,...,x_n): = \sum ^n a_{ij}x_ix_j+
\sum ^n b_ix_i +c=0\}$, where $a_{ij}=a_{ji}$. Then for any $\epsilon _1,
..., \epsilon _{n-1} \in \{ -1, 1\}$ we have $0=(F(\epsilon _1 d_1, ...,
\epsilon _{n-1} d_{n-1}, d_n)-F(\epsilon _1 d_1,..., \epsilon _{n-1} d_{n-1},
-d_n))/2=d_n(\sum ^{n-1} 2\epsilon _i a_{in}d_i +b_n)$, hence $0=\sum ^{n-1}
2\epsilon _i a_{in} d_i +b_n$, that readily implies 
$a_{in}=b_n=0\,\,(i\le n-1)$. 
Similarly one shows $a_{ij} =0 \,\,(i<j \le n)$ and 
$b_i =0 \,\,(i \le n)$, that implies the statement of the lemma. 
\end{proof}

\begin{corollary} Let $E \subset 
{\R} ^n$ be an ellipsoid with axes of 
different lengths, and $B \subset {\R} ^n$ a box of size $a_1 \times ...
\times a_n$, where $a_1=...=a_{n_1}<a_{n_1+1}=...=a_{n_1+n_2}<...<a_{n_1+
...+n_{k-1}+1}=...=a_n$ and $n=n_1+...+n_k$. Then there are $n!/(n_1!...n_k!)$
ways of inscribing a similar copy $B'$ of $B$ into $E$ (and in each of these 
cases the centre of $B'$ coincides with that of $E$).
\label{ellipsoid}
\end{corollary}

\begin{proof} By Lemma~\ref{center} 
the centre of an inscribed similar copy $B'$ of
$B$ is at the centre of $E$, and its edges of length $a_{n_1+...+n_{i-1}+1}$ 
are 
parallel to some $n_i$ axes of $E$. This gives a partition of the $n$ axes of 
$E$ to classes of sizes $n_1,...,n_k$; the number of such partitions is 
$n!/(n_1!...n_k!)$. Fixing one such partition, $B'$ is determined, up to 
magnification from its centre. Its size is determined by the requirement that 
a given vertex of it should belong to $\partial E$ (and then, by symmetry, all 
of its vertices belong to $\partial E$). 
\end{proof}

\begin{lemma} 
\label{transversal}
Let $E \subset {\R} ^n$ be an ellipsoid with axes of 
different lengths, and of equation $F(x_1,...,x_n):=\sum ^n a_{ii} x_i ^2
-1=0$. Let $P$ be a box inscribed to $E$, of vertex set $\{v_1,..., v_{2^n}\}=
\{(\pm d_1,...,\pm d_n)\}$. Let us consider the analytic map $f$, sending $A 
\in SO(n)$ to $(F(Av_1),...,F(Av_{2^n})) \in {\R} ^ {2^n}$. Then at $A=I\,\,\,
f$ has full rank $n(n-1)/2$, and the diagonal $\Delta$ of ${\R} ^ {2^n}$ 
intersects the tangent space of $f(SO(n))$ at $f(A)=f(I)$ transversally.
\end{lemma}

\begin{proof} 
We have $f(I)=0$, and $T_I(SO(n))$ is the set of skew symmetric 
$n \times n$ matrices. Let $(d\alpha _{ij}) \in T_I(SO(n))$. Then $f(I+
(d\alpha _{ij}))=(F(v_1+(d\alpha _{ij})v_1,...,F(v_{2^n}+(d\alpha _{ij})
v_{2^n}))$. We have to show that for $0 \ne (d\alpha _{ij}) \in T_I(SO(n))
$ we have that $f(I+(d\alpha _{ij}))$ does not lie on the diagonal of $T(
{\R} ^ {2^n}) \cong {\R} ^{2^n}$ (thus, in particular, is not $0$). We will 
argue indirectly.
 
For $\{v_1,...,v_{2^n}\} \ni v =(x_1,...,x_n)$ we have $F(v+(d\alpha _{ij})v)=
F(v)+F' \cdot ((d\alpha _{ij})v)=2\sum _i a_{ii} x_i(\sum _j d\alpha _{ij} 
x_j)=2\sum _{i<j}(a_{ii} -a_{jj})d\alpha _{ij}x_ix_j$, taking in account the 
skew symmetry of $(d\alpha _{ij})$. Here by hypothesis $a_{ii} -a_{jj} \ne 0$ 
for $i<j$. If $F(v_1+(d\alpha _{ij})v_1)=...=F(v_{2^n}+(d\alpha _{ij})
v_{2^n})$, then the quadratic form $\sum _{i<j} (a_{ii} -a_{jj})d\alpha _{ij}
x_ix_j$ 
assumes equal values for each $(x_1,...,x_n)=(\pm d_1,...,\pm d_n)$. Like in 
Lemma~\ref{center}, 
this implies $(a_{ii}-a_{jj})d\alpha _{ij}=0$ for each $i<j$, hence 
$(d\alpha _{ij}) =0$ as asserted. \end{proof}

\paragraph{\it Proof of Theorem~\ref{cube}.} 
For any even map $F:S^2\rightarrow \R$
we define the map $f:SO(3) \to {\R}^4$ sending
$A \in SO(3)$ to $(F(\pm Av_1), ...,F(\pm Av_4))$. We are going to show 
that the image of $f$ intersects the diagonal $\Delta$ of ${\R}^4$,
which is equivalent to the statement of Theorem~\ref{cube}. 

Of course maps 
from $SO(3)$ to $\R^4$ do 
exist whose images do not intersect the diagonal. Therefore
we need some extra condition on the map $f$. Namely $f$ as defined above
should respect the rotation group of the cube in the sense that it
should be equivariant with respect to the free
$S_4$ action $\rho_{S_4}$ on $SO(3)$ 
(given by right multiplication by $S_4\subset SO(3)$) and 
$\tilde{\tau}_{S_4}$ on $\R^4$ (given by permuting the coordinates), i.e., 
$\tilde{\tau}_{S_4} (g)f=f\rho_{S_4} (g)$ for each $g \in S_4$. 
We show that the image of 
such an equivariant map intersects the diagonal $\Delta$ 
of $\R^4$. 

Indirectly suppose that an equivariant map 
$g:(SO(3),\rho_{S_4})\rightarrow (\R^4,\tilde{\tau}_{S_4})$ exists, whose
image
does not intersect the diagonal. Applying the projection along $\Delta$ to
$\Delta^\perp\cong \R^3$ (the orthogonal complement of $\Delta$) of $\R^4$ 
we get an equivariant map 
$p:(\R^4,\tilde{\tau}_{S_4})\rightarrow (\R^3, \tau_{S_4})$, where
the $S_4$ action $\tau_{S_4}:S_4 \to O(3)$ is 
the action of the symmetry group of 
the regular tetrahedron in $\R^3$, the four letters corresponding to the four 
vertices of the regular tetrahedron. 
 Thus we have an equivariant map
$p\circ g: (SO(3),\rho_{S_4})\rightarrow (\R^3,\tau_{S_4})$ whose image
does not contain the origin. Finally we have the equivariant map 
$$h:(\R^3\setminus \{0\}, \tau_{S_4})\rightarrow (S^2,\tau_{S_4}),$$ which
is the projection along the 
radial direction. Therefore we have an equivariant map
$h\circ p \circ g:(SO(3),\rho_{S_4})\rightarrow (S^2,\tau_{S_4})$. 
The following proof of Proposition~\ref{s4} rules out the
possibility of existence of such mappings, completing the proof of
Theorem~\ref{cube}. $\square$

\paragraph{\it Proof of Proposition~\ref{s4}.} 
An equivariant map $(SO(3),\rho_{S_4})\to (S^2,\tau_{S_4})$
would give by construction a section of the
$S^2$  bundle $\eta_{\tau_{S_4}}$ defined as the sphere bundle 
$$S^2\rightarrow (SO(3)\times S^2)/(\rho_{S_4}\times \tau_{S_4})
\rightarrow SO(3)/\rho_{S_4}.$$ We show that this sphere bundle does
not have a section by showing that
the third Stiefel-Whitney class of the sphere bundle
$w_3(\eta_{\tau_{S_4}})\neq 0$  
does not vanish. For this
it is sufficient to show a section of the corresponding vector bundle 
$$\bar{\eta}_{\tau_{S_4}}:=(SO(3)\times \R^3)/(\rho_{S_4}\times \tau_{S_4})
\rightarrow SO(3)/\rho_{S_4}$$ 
which intersects transversally 
the zero section in an odd number of points. Namely it is well known
that the 
cohomology class of the zero set of a section which is transversal to
the zero section coincides with the Euler class of the vector bundle
which in our case is the third Stiefel-Whitney class $w_3$. 

For exhibiting such a section we consider an ellipsoid $E\subset \R^3$,
centred at the origin, with axes of different lengths. 
Moreover we let 
$s_E:(SO(3),\rho_{S_4})\rightarrow (\R^3, \tau_{S_4})$ be the section of
$\bar{\eta}_{\tau_{S_4}}$ corresponding to the convex body $E$. 
 By Corollary~\ref{ellipsoid} $E$ admits exactly one inscribed cube.
Thus $s_E$ vanishes at exactly one point of $SO(3)/\rho_{S_4}$, 
namely at the point corresponding to this 
inscribed cube. By Lemma~\ref{transversal}, at this point the intersection is
transversal. 
   
Thus $w_3(\eta_{\tau_{S_4}})\neq 0$ indeed,
completing the proof of Proposition~\ref{s4}, 
hence of Theorem~\ref{cube}, and of Corollary~\ref{inscribe}. $\square$


\section{$G$-equivariant maps} 
\label{equivariant}

The theorem above gives that there is no $S_4$-equivariant map from
$SO(3)$ to $S^2$. One may wonder for which subgroup of $S_4$ we can find
an equivariant map. To be more precise we need some more notation. 

A subgroup $G\subset S_4$ gives the $G$-actions $\rho_{G}$ and $\tau_{G}$
by restricting the $S_4$ actions $\rho_{S_4}$ and $\tau_{S_4}$ to the 
subgroup $G$. Moreover we get the $S^2$ bundle $\eta_{\tau_G}$ as the
sphere bundle $S^2\rightarrow (SO(3)\times S^2)/(\rho_G\times \tau_G) 
\rightarrow SO(3)/\rho_G$. Our problem is thus to determine for which 
subgroup $G$ of $S_4$ we get a section of the $S^2$ bundle 
$\eta_{\tau_G}$.       

The \cite{gap} algebraic program 
package tells that, up to conjugation, there are 11
different subgroups of $S_4$; these are given by it by generators, and are 
isomorphic to the cyclic groups $C_1,\,\,C_2,\,\,
C_2,\,\,C_3,\,\,C_4$, the dihedral
groups $D_2,\,\,D_2,\,\,D_3,\,\,D_4$, the alternative group $A_4$, and $S_4$,
respectively. Since $A_4$ occurs in this list once, and is a normal subgroup 
of $S_4$, there is only one subgroup of $S_4$, isomorphic to $A_4$, and, for 
any subgroup $G$ of $S_4$, the property $G \subset A_4$ is invariant under 
taking conjugates of $G$. Thus, by the above list, the subgroups of $S_4$, not
contained in $A_4$, are, up to conjugation, the following: $[(ij)] \cong
C_2$, $[(ijkl)] \cong C_4$
(corresponding to the group of rotations through multiples of
$\pi /2$ about a coordinate 
axis), $[(ij),(kl)] \cong D_2 \cong C_2 \times C_2$, $D_3 \cong S_3$ the 
subgroup of permutations fixing some $i$, $D_4$ the rotation group of a 
non-cubical square-based box (when $S_4$ is represented in $SO(3)$ as above),
and $S_4$. ($[B]$ is the subgroup generated by $B \subset S_4$, and $i,j,k,l$
are different elements of $\{1,...,4\}$.)

\vskip.4cm

\noindent {\it 1. Case $G\cong S_4$} 

\vskip.4cm

The previous theorem has showed that in this case 
$\eta_{\tau_{S_4}}$ does not have a section.


\subsection{Inscribed square based boxes}
\label{squares}
    
\noindent {\it 2. Case $G\cong D_4$} 

\vskip.4cm
In this case we do not have a section as the following proposition
shows. Still we note that the following three statements are 
generalizations of Proposition~\ref{s4}, 
Theorem~\ref{cube} and Corollary~\ref{inscribe}, 
respectively. Like Theorem~\ref{cube},
also Theorem~\ref{square} follows from Theorem D of \cite{griffiths}, too.

\begin{proposition}
\label{d4} 
There is no equivariant map from $(SO(3),\rho_{D_4})$ to 
$(S^2,\tau_{D_4})$.
\end{proposition}

\begin{proof} The proof of this proposition follows the proof 
of Proposition~\ref{s4},  with the exception that we need the following
consequence of Corollary~\ref{ellipsoid}:
there are three ways of inscribing a non-cubical square based box
with given ratio of edges into an ellipsoid in ${\R}^3$ with  axes of 
different lengths. 

Now the proof of the proposition is complete noting that $3$ is odd and 
therefore $w_3(\eta_{\tau_{D_4}})$ does not vanish.  
\end{proof}

\begin{theorem} 
\label{square}
Let $F :S^2 \to {\R}$ be a continuous even map, 
and let $\{
\pm v_i \mid 1 \le i \le 4\}$ be the vertex set of a square based box inscribed
to $S^2$. Then there exists $A \in SO(3)$ such that $F(\pm Av_1)=F(\pm Av_2)
=F(\pm Av_3)=F(\pm Av_4)$.
\end{theorem}

\begin{proof} 
As in the case of the cube we can define a map $f:SO(3) \to
{\R}^4$ by setting $f(A)=(F(\pm Av_1), ..., 
F(\pm Av_4))$. We have to show that the image of $f$ intersects the diagonal 
$\Delta$ in ${\R}^4$. Notice that the rotation group of a square based box
contains $D_4\,\, (\subset S_4$) 
(and in general equals $D_4$). Now an identical
argument as in the case of the cube in Theorem~\ref{cube} 
yields the statement of 
this theorem, using Proposition~\ref{d4} instead of Proposition~\ref{s4}. 
\end{proof}

\begin{corollary}
\label{inscribesquare}
Every centrally symmetric convex body $K$ in $\R^3$ admits an inscribed
square based box, with any given ratio of the height to the basic edge, 
and centred at the 
centre of $K$.
\end{corollary}

\begin{proof} Let us consider a square based box, with the given ratio of the 
height to the basic edge, 
that is inscribed to $S^2$, and has vertex set $\{ \pm v_i \mid
1 \le i \le 4\}$.
Like at the reduction of Corollary~\ref{inscribe} to Theorem~\ref{cube}, 
we suppose that the centre of $K$ is the 
origin, and we let $\|\,\,\|_K$ be the Minkowski norm in ${\R}^3$, 
associated to $K$. We apply Theorem~\ref{square} 
for $F=\|\,\,\|_K \mid_{S^2}$ and the $v_i$ 
chosen above, obtaining $A \in SO(3)$ such that $(\| \pm Av_1 \| _K, 
..., \| \pm Av_4 \|_K) $ is a point
$(\lambda, ..., \lambda)\,\,(\lambda >0)$  
on the diagonal $\Delta$ of ${\R}^4$. Then 
the square based box with vertex set $\{\lambda ^{-1} (\pm Av_i) \mid
1 \le i \le 4\}$ is inscribed in $K$, has the given ratio of the height to
the basic edge, 
and is centred at the centre of $K$. 
\end{proof}


\subsection{Existence of $G$-equivariant maps}

\label{existence}

\noindent {\it 3. Case $G \subset A_4$}

\begin{theorem}
For $G \subset A_4$ there exists a $G$-equivariant
map from $(SO(3), \rho _G)$ to $(S^2, \tau _G)$.
\end{theorem}

\begin{proof} 
We will show
that the $S^2$-bundle
$\eta_{\tau_{A_4}}$ defined at the beginning of Section~\ref{equivariant}
is actually trivial. Then the $S^2$-bundle $\eta _{\tau _{A_4}}$ will have a 
section, that implies by construction the existence of  
$A_4$-equivariant maps, that are also $G$-equivariant for any $G \subset A_4$.

Think of 
$A\in SO(3)$ as an orthonormal basis. Let $p_i:SO(3)\rightarrow S^2$ be
the forgetful map associating to $A\in SO(3)$ its $i$'th basis vector. 
These maps are $A_4$-equivariant as $A_4$ is the rotation group of the 
regular tetrahedron. In this way we get three linearly independent
sections $p_1, p_2$ and $p_3$ of the vector bundle 
$\bar{\eta}_{\tau_{A_4}}$. Thus $\eta_{\tau_{A_4}}$ is trivial indeed.  
(Another way of seeing this is the following. If $G$ is a discrete subgroup 
of $SO(3)$, then the bundle $\bar{\eta}_{\tau_G}=
SO(3) \times _G \R ^3=T(SO(3)/G)$ 
is the tangent bundle of $SO(3) / G$ (left cosets), and this is 
 trivial. Namely the tangent bundle of the quotient of a Lie group with 
respect to a discrete subgroup is trivial.)
\end{proof}

\begin{remark} This argument does not work for $G \not\subset A_4$,
because in 
these cases the action of $G$ is given by $g(\hat g,x)=(g \hat g,\pm  gx)$, where the sign is positive for $g \in A_4$ and negative 
for $g \in S_4 \setminus A_4$.  
\end{remark}

In what follows, we will investigate subgroups $G \not\subset A_4$. 

\vskip.4cm

\noindent {\it 4. Case $G \cong C_2,D_2,D_3$}

\vskip.4cm 

In the following theorem we will use the notations introduced at the 
beginning of Section~\ref{equivariant}. 

\begin{theorem}
Let $G=[(ij)] \cong C_2,\,\,G=[(ij),(kl)] \cong D_2$
or let $G \,\,( \cong 
D_3 \cong S_3)$ be the subgroup of permutations fixing some $i$. Then there
exists a $G$-equivariant map from $(SO(3), \rho _G)$ to $(S^2, 
\tau_G)$.
\end{theorem}

\begin{proof} 
Recall that $\tau _G :G \to O(3)$ is the action of the subgroup 
$G$ of the symmetric group $S_4$ of the regular tetrahedron, the four letters 
corresponding to the four vertices of the regular tetrahedron. In other words,
for $g \in G$ we have that $\tau _G (g)$ is that element of $O(3)$, whose 
restriction to the vertex set of the regular tetrahedron inscribed to $S^2$ 
equals the permutation $g \in G$ of these vertices. 

For $G=[(ij)]$ any $x \in S^2$, equidistant to the $i$'th and $j$'th vertices,
is a fixed point for $\tau _G$. For $G=[(ij),(kl)]$ any vector $x \in S^2$, 
orthogonal to the lines connecting the $i$'th and $j$'th vertices, and the
$k$'th and $l$'th vertices of the tetrahedron, respectively, is a fixed point 
of $\tau _G$. For $G$ being the subgroup of permutations fixing some $i$, 
the $i$'th vertex is a fixed point of $\tau _G$. Then a $G$-equivariant map can
be given in all three cases as a constant map $(SO(3), \rho _G) \to (S^2,
\tau _G)$, having as value a fixed point $x$ of $\tau _G$. 
\end{proof}

\noindent{\it 5. Case $G \cong C_4$}
 
\begin{theorem} There does exist a $C_4$-equivariant 
map from $(SO(3),\rho_{C_4})$ 
to $(S^2,\tau_{C_4})$.
\end{theorem}

\begin{proof} 
We show that there is a section of the $S^2$ bundle 
$\eta_\tau=(SO(3)\times S^2)/(\rho_{C_4}\times \tau_{C_4})$ over
the orbit space $L_8=SO(3)/\rho_{C_4}$. This is equivalent to 
the statement of the theorem. 

Obstruction theory (cf. \cite{milnor-stasheff}, 
pp. 139-143) 
tells us that the existence of such a section is equivalent
to the vanishing of the obstruction class 
$o_3(\eta_\tau)\in H^3(L_8,\pi_2 (S^2))$ 
sitting in a cohomology space with twisted  coefficients. We show the
vanishing of $o_3(\eta_\tau)$. We use a theorem of Steenrod 
(cf. \cite{steenrod}, 38.8) which claims that there is a homomorphism 
$$\delta^*:H^2(L_8,\pi_1(SO(3))\rightarrow H^3(L_8,\pi_2(S^2))$$ which sends
the second obstruction class $o_2(\eta_\tau)$ to $o_3(\eta_\tau)$. 
Therefore it is sufficient to
show that the second obstruction class vanishes. On the other hand
the second obstruction class $o_2(\eta_\tau)$ coincides 
with the second
Stiefel-Whitney class of $\eta_\tau$, i.e., $o_2(\eta_\tau)=w_2(\eta_\tau)$
 (cf. \cite{milnor-stasheff}, p.140). 
Thus the proof of $w_2(\eta_\tau)=0$ gives the theorem. 

The vector bundle $\eta_\tau=\eta_{\tau_1}\oplus \eta_{\tau_2}$ 
decomposes as the direct sum
of a rank $1$ and a rank $2$ vector bundle, where $\tau_1$ is the non-trivial
representation of $C_4$ on $\R^1$ and 
$\tau_2$ is an effective representation 
of $C_4$ on $\R^2$ (the generator is sent to the rotation through the angle 
$\frac{\pi}{2}$), and 
$\eta_{\tau_i}=(SO(3)\times \R^i) / (\rho\times \tau_i)$. 
Thus we get that 
$$w_2(\eta_\tau)=w_1(\eta_{\tau_1})\cup w_1(\eta_{\tau_2}) + 
w_2(\eta_{\tau_2}).$$

This formula shows that the vanishing of $w_1(\eta_{\tau_2})$ and 
$w_2(\eta_{\tau_2})$ yields the vanishing of $w_2(\eta_{\tau})$. 
The first Stiefel-Whitney class
$w_1(\eta_{\tau_2})$ vanishes since 
the representation $\tau_2$ is orientation preserving
and hence $\eta_{\tau_2}$ is oriented. 
The second Stiefel-Whitney class $w_2(\eta_{\tau_2})$ vanishes since, 
as we are going to show, 
$$\eta_{\tau_2}=\eta_\gamma\otimes\eta_\gamma$$ with some 
$S^1$ principal bundle 
$\eta_\gamma$ (here $\eta_{\tau_2}$ is considered as 
an $S^1$ 
principal bundle).
Hence the Euler class 
$e(\eta_{\tau_2})=2e(\eta_\gamma)$ is even, yielding the desired 
vanishing of $w_2(\eta_{\tau_2})$.

Thus we are left with constructing an $S^1$ 
principal bundle $\eta_\gamma$ with 
$\eta_{\tau_2}=\eta_\gamma\otimes \eta_\gamma$. The universal
covering of $L_8$ is the composite covering 
$S^3\rightarrow \R P^3\cong SO(3)\rightarrow L_8$. Moreover, we claim that this 
covering is given by the $C_8$ action $\tilde{\rho}$ 
on $S^3$, sending the generator $\alpha$
to the transformation 
$q\mapsto e^{\frac{2\pi i}{8}}q$ of $S^3$, the unit sphere
of the quaternionic algebra $\Q$. To see this, divide out first by the 
$C_2\subset C_8$ action defined by $\alpha^4(q)=e^{\frac{2\pi i}{2}}q=-q$. The
quotient is 
$\R P^3$ which can be identified with $SO(3)$, by associating to any 
element $(q,-q)\in \R P^3$ a special orthogonal transformation of 
${\rm Im}(\Q)=\R^3$, namely the conjugation with $q$. Now we see that 
the $C_4=C_8/C_2$ 
action on $\R P^3$ inherited from the $C_8$ action $\tilde{\rho}$
of $S^3$ corresponds to the 
$C_4$ action $\rho$ of $SO(3)$ with respect to the axis 
$i\in {\rm Im}(\Q)$ (i.e., to the rotation of the space Im$(\Q) \cong
{\R} ^3$ through $\pi /2$ about the axis $i$),
via the above 
identification. Hence indeed $L_8$ is a lens space. 
    
Let $\eta_{\gamma}=(S^3\times \R^2)/(\tilde{\rho}\times \gamma)$ 
be the principal $S^1$ bundle on $L_8$ given by the 
representation $\gamma$ of $C_8$ on $\R^2$, sending the generator to the 
rotation through the angle $\frac{\pi}{4}$. 
Now the relation $\eta_\gamma\otimes\eta_\gamma= \eta_{\tau_2}$ is 
immediate. The proof is complete.   
\end{proof}


\section{Inscribed boxes}
\label{box}

Now we will prove a theorem that contains both our previous 
Corollary~\ref{inscribe}
and its more general form Corollary~\ref{inscribesquare} 
as special cases. However, as 
mentioned above, its proof uses the quite involved Theorem D in 
\cite{griffiths},  
rather than its special cases Theorem~\ref{cube}, and its more general form 
Theorem~\ref{square}.

\begin{theorem}
\label{boxes}
Every centrally symmetric convex body $K$ in ${\R}^3$ admits an inscribed 
box, similar to any given box, and centred at the centre of 
$K$.
\end{theorem}

\begin{proof} Let us consider a box similar to the given box,
that is inscribed to $S^2$, and one of whose faces has vertex set 
$\{v_1,..., v_4\}$.
Like at  Corollaries~\ref{inscribe} and \ref{inscribesquare}, 
we suppose that the centre of $K$ 
is the origin, and we let $\|\,\,\|_K$ be the Minkowski norm in ${\R}^3$,
associated to $K$. We let $F=\|\,\,\|_K \mid_{S^2} :S^2 \to 
{\R}$. Then Theorem D of \cite{griffiths}
(quoted in Section~\ref{introduction}) 
implies for $F:S^2 \to {\R}$ and the set 
$\{v_1,...,v_4\}$ 
of the vertices of a planar rectangle inscribed to $S^2$, that 
there exists $A \in SO(3)$ such that $((\|Av_1\|_K,..., \|Av_4\|_K)=
\nolinebreak )\,\,
(F(Av_1),...,F(Av_4))$ is a point $(\lambda,...,\lambda)\,\,(\lambda >0)$ 
on the diagonal $
\Delta$ of ${\R}^4$. Then the box with vertex set $\{\lambda ^{-1}(\pm
Av_i) \mid 1 \le i \le 4\}$ is inscribed in $K$, is similar to the given box, 
and is centred at the centre of $K$. 
\end{proof}


\section{Inscribed cubes in $C^1$ surfaces}
\label{konkavfejezet}

In this section we prove theorems about inscribability of cubes and square
based boxes to certain centrally symmetric surfaces, that do not follow from
Theorem D of \cite{griffiths}. We will roughly follow the proofs of Theorems
1 and 2 of \cite{jerrard}, showing that a plane analytic Jordan curve admits 
an inscribed square, and in the "general" case (precised there) their number 
is odd. Cf. also \cite{gromov}, 
showing inscribability of homothets of simplices 
to hypersurfaces in ${\R}^n$, under topological and smoothness assumptions,
and also solving a Knaster type problem.
Another related theorem is \cite{griffiths}, Theorem C, stating that each 
$C^2$--embedded $S^2$ in $\R ^3$ admits for each $\rho >0$ an inscribed "skew 
box of aspect $\rho $". This is obtained from a square based right prism with 
ratio of height to basic edge equal to $\rho $, by rotation of one base in 
its own plane about its centre through some angle.

First we prove Theorem~\ref{konkav}:

\paragraph{\it Proof of Theorem~\ref{konkav}.}
By Corollary~\ref{ellipsoid} the statement is true if $f(S^2)$ is an
ellipsoid with pairwise different lengths of axes. Moreover there is a single
cube on such an ellipsoid and this is ''transversally true", i.e., remains true
for any sufficiently small perturbation 
in the $C^1$ topology
(invariant under the map $x \mapsto
-x$) of the ellipsoid (cf. Lemma~\ref{transversal}).

\paragraph{The plan of the proof.} 
We want to reduce the general statement of the
theorem to the above special case. Namely we show that the cubes on the surface
$f(S^2)$ are in bijective correspondence with the intersection points of two
cycles ($\alpha$ and $\beta _f$) 
of complementary dimensions of a $C^1$ manifold
with boundary. Here
$\alpha$ is formed by (the equivalence classes of) those eight-tuples of points
in $\R^3 \setminus \{0\}$ 
which are vertices of a cube (centred at $0$). The cycle $\beta _f$ is
formed by (the equivalence classes of) those eight-tuples of points
which are all on the surface $f(S^2).$
We may assume that their intersection is transversal, and consists of
finitely many points, since else it is certainly not empty, and then there is 
nothing to prove.
Now if we replace the surface $f(S^2)$ by any other one -- $g(S^2)$, say, --
satisfying the same assumptions and thus representing
the same homotopy class, 
then the cycle $\beta _f$
will be replaced by a cycle $\beta _g$, homological to it. Since $\alpha$ is a 
cycle,
the intersection number (mod $2$) $\alpha \cap \beta _g$ remains equal to
$\alpha \cap \beta _f.$ 
On the other hand, $\alpha \cap \beta _g$ is the number (mod $2$)
of the cubes on $g(S^2).$ If we choose 
$g(S^2)
= $ an ellipsoid with different lengths of axes, with $g(-x)=-g(x)$ 
(let $g$ be the radial projection of $S^2$ to the surface of the ellipsoid)
then, 
by the special case
proved earlier, $\alpha \cap \beta _g \ne 0$ mod $2.$

\paragraph{Preliminary remarks.}
The (oversimplified) plan described above will be modified by the following
ways and following reasons.
\vskip.5cm
1. Since $f(S^2)$ is invariant under the multiplication by $-1$ it will be
sufficient to find four vertices $v_1, v_2, v_3, v_4$ of a face of a cube,
in this cyclic order,
because then the eight vectors $\pm v_i$ automatically will
form a cube on $f(S^2)$.
(In this proof by a cube 
we always mean one centred at the origin.) Therefore we will form
the cycle $\alpha$ not from the $8$-tuples of points but from quadruples,
forming the vertices of a face of a cube. Similarly $\beta _f$ 
will be formed not
from $8$-tuples but quadruples of points on $f(S^2)$. The dimensions of these 
cycles in $({\R}^3)^4$ are $4$ and $8$, respectively, thus are complementary.
\vskip.5cm
2. Naturally we do not consider two cubes  different, if they  differ only
by a permutation
of the vertices. (Otherwise we would get each cube $48$
times
and therefore the homological intersection index with $\Z_2$ coefficients would
not give any information on their number.) Therefore we have to factorize
out $(\R^3)^4$ by the "automorphism group of the cube", i.e., by the group 
$G:=\{ T \in O(({\R}^3)^4) \mid (Tv_1, Tv_2, Tv_3, Tv_4) $ is a permutation of 
$(v_1, v_2, v_3, v_4)$, or of $(v_i, v_{i+1}, -v_{i+3}, -v_{i+2}),\,\,
1 \le i \le 4$, indices meant cyclically, or of $(-v_1, -v_2, -v_3, -v_4)$, in 
each of these cases 
either {\it {preserving}} or {\it {reversing}} 
their above given cyclic order$\}$.
(Observe that if $v_1, v_2, v_3, v_4 \in {\R}^3 $ 
are the vertices of a face of 
a cube, in this cyclic order, then in a natural way these transformations
bijectively correspond to the elements of the symmetry group ($ \subset 
O({\R}^3) $) of this cube.)
This raises the problem
that this action is not free at certain points  
of the set $\{(v_1,v_2,v_3,v_4) \in ({\R ^3})^4 \mid v_i=\pm v_j {\rm {\,\,for
\,\,some\,\,}} 1 \le i < j \le 4 \}$. Since later we will have to exclude
quadruples $(v_1,...,v_4) \in ({\R ^3})^4$ having a zero coordinate $v_i$,
we define the 
{\it {singular set}} $S$ as the larger set 
$$S:=\{(v_1, v_2, v_3, v_4) \in ({\R}^3)^4 \mid
v_i = 0 {\rm {\,\, for\,\, some\,\,}} 1 \le i \le 4 $$
$${\rm {or\,\,}} 
v_i = \pm v_j {\rm {\,\, for\,\, some\,\,}} 1 \le i <
j \le 4\} \subset (\R^3)^4.$$
The action of $G$ is free on $({\R }^3)^4 \setminus S$. 
By non-freeness of the action on $({\R ^3})^4$ 
the quotient space $({\R ^3})^4/G$ will not be a manifold,
only an orbifold.
\vskip.5cm
Let $A$ be the set corresponding to the 
cycle $\alpha$, i.e.,
\begin{eqnarray*}
A := \{(v_1, v_2, v_3, v_4) 
\mid v_i \in \R^3 \setminus \{0\} \ \mbox{\rm and}\ v_1, v_2, v_3, v_4
\ \mbox{\rm are the vertices of a face of a cube}, \\ 
\mbox{in this cyclic order}\}.
\end{eqnarray*}
Note that $A$ is a cone in $(\R^3)^4$ with its 
vertex at the origin, with its vertex deleted, and with no
point on the
set $S$.

Of course $X:=({\R}^3)^4 \setminus S \subset ({\R}^3)^4$ is an open set, 
invariant under $G$, 
on which $G$ acts freely. Thus we can consider the quotient manifold 
$
X/ G$. 
But the trouble is that this manifold is 
not compact.

\vskip.5cm
3. So the final modification is that from the very beginning we shall work
not with the whole set $X$, but only with some sufficiently large 
compact subset of it
(to be specified later).
The set $X$ is an open set invariant under $G$, and $G$ acts freely on it. Let 
$\varepsilon _1 >0$ be a sufficiently small number, $K>0$ be sufficiently large
(to be specified later), and let $\varepsilon _2 = \varepsilon _1 /{100}$. Let
$\varphi : X \to \R,$
$$\varphi (v_1,...,v_4):= {\rm {max}} \left \{ \frac
{\varepsilon _1 }{\|v_i \|}\,\,\,\,(1 \le i \le 4), \,\,\frac{\|v_i\|}{K}\,\,\,\,
(1 \le i \le 4),\,\,\frac{\sqrt {2} \varepsilon _2 }{\|v_i \pm v_j\|}\,\,\,\,
(1 \le i < j \le 4) \right \}.$$
Observe that $\|v_i \pm v_j \|/\sqrt {2} $ is the distance of $(v_1,...,v_4)$ 
to the subspace given by the equation $v_i \pm v_j =0$ (since $(v_1,...,v_4)
\mapsto \left( (v_1 +v_2 )/\sqrt {2} , (v_1 -v_2 )/\sqrt {2}, v_3, v_4 \right) $
is an orthogonal transformation). Then $\varphi $ is continuous, is invariant 
under $G$, and for any $r>0$ the set $\varphi ^{-1} ([0,r])$ is a compact 
subset of $X$. However, we do not know whether $ \varphi ^{-1} ([0,r])
\setminus {\rm {int}}\varphi ^{-1}([0,r])$
is a topological manifold for a generic $r>0$. Therefore we define
$\psi : X \to \R$,
$$\psi (v_1,...,v_4):= \left( \sum _{i=1} ^4 \frac
{\varepsilon _1 ^2}{\|v_i\|^2}+\sum _{i=1} ^4 \frac {\|v_i\|^2}{K^2} +
\sum _{1 \le i <j \le 4} \frac {2 \varepsilon _2 ^2}{\|v_i \pm v_j\|^2}
\right) ^{1/2} .$$
Then $\psi $ is an analytic, hence $ C ^ {\infty } $ function on $X$, is 
invariant
under $G$, and
$$ \varphi \le \psi \le \sqrt {20} \varphi .$$

By Sard's lemma, almost all $r \in \R $ is not a critical value of $ \psi $
(i.e.,
for each $(v_1,...,v_4) \in \psi ^{-1} (r )$ we have $\psi ' (v_1,...,
v_4) \ne 0$ -- that means 
at least one partial derivative is different from $0$ --, 
and hence $\psi ^{-1} (r )$ is a closed embedded 
$C^{\infty}$ $11$-submani\-fold of $X$, 
possibly empty), cf. \cite{howard},
p. 66. Let $r \in (\sqrt {20},5)$ be a not critical value of $\psi $, and let
$$X_r :=  \psi ^{-1} \left( [0,r )\right) = 
\psi ^{-1} \left( (0,r )\right) \,\,(\supset 
\varphi ^{-1} \left( [0,1 ] \right) ). 
$$
The set $Y_r:=\psi ^{-1} (r ) \subset \varphi ^{-1} \left( (1,5 ) \right)
$ is a non-empty closed 
embedded $C^{\infty }$ $11$-submanifold of $X$, and is invariant under $G$. 
Further, $Y_r = \partial X_r$, since $r $ is not a 
critical value of $\psi $. 

The factor space $\overline{X_r} /G = (X_r \cup Y _r
)/G$ is a compact $C^{\infty }$ manifold with boundary $ 
Y_r/G$. Since both of the sets $A,
f(S^2)^4 \subset ({\R}^3)^4$ are $G$-invariant $C^1$ manifolds, so their 
intersections with $X_r$ and $Y_r$ are $G$-invariant as well. 
The sets $A  \cap X_r,\,\,f(S^2)^4 \cap 
X_r$ are $C^1$ manifolds as well (but we do not know if $A
 \cap {\overline{X_r}},\,\,f(S^2)^4 \cap 
{\overline{X_r}}$ are $C^1$ manifolds with boundary, since we do not know
anything about the sets $A \cap Y_r,\,\,f(S^2)^4 \cap Y_r)$. 
Thus we can form the 
quotient $C^1$ manifolds $[A  \cap X_r]/G,
\,\,[f(S^2)^4 \cap X_r]/G$ (observe that $A \cap X_r,\,\,f(S^2)^4 \cap X_r$
are disjoint to the singular set $S$), and the quotient sets $$A_r:=
[A \cap {\overline{X_r}}]/G,\,\,B_r
:=[f(S^2)^4 \cap {\overline{X_r}}]/G $$ (about which we do not know if 
they are $C^1$ manifolds with boundary). For the statement in the last brackets
it would be necessary to know if $Y_r /G=\partial (X_r/G)$ 
intersects the manifold $A/G$ and the orbifold $f(S^2)^4/G$, that is a 
manifold in a neighbourhood of $Y_r /G$,
transversally. 

However, by Thom's transversality theorem (cf. \cite{hirsch}
Ch. 3, Theorem 2.5) this can be achieved, if we replace $Y_r/G$ by 
a $C^1$ manifold $Q_r \subset X / G$, that is a small, generic
$C^1$ perturbation of $Y_r/G$.
Moreover, we may assume that
$$Q_r={\overline {P_r}} \setminus P_r = \partial {\overline {P_r}},$$
for some open set $P_r \subset X/G$, whose symmetric 
difference with $X_r/G$ lies in a small tubular neighbourhood of 
$Y_r/G$, included in $[\varphi ^{-1} ((1, 5))]/G$.
Further, we have 
$$Q_r \subset \varphi ^{-1}\left( [1,5] \right) /G \,\, {\rm and}\,\,
P_r \subset \varphi ^{-1}\left( [0,5)\right) /G.$$
Moreover, the modified sets 
$$A'_r:=[A/G] \cap 
{\overline{P_r}} \,\,{\rm {and}}\,\, B'_r:=[f(S^2)^4/G] \cap {\overline{P_r}}
=[\left( f(S^2)^4 \setminus S \right) /G] \cap {\overline{P_r}}$$
are $C^1$ manifolds 
with boundary, their boundaries $[A/G] \cap 
{{Q_r}}$ and $[\left( f(S^2)^4 \setminus S \right) /G] 
\cap {{Q_r}}$
lying in $\partial {\overline {P_r}}$; in other words,
they are relative cycles in the chain complex 
${\cal C} _* ({\overline {P_r}}, \partial {\overline{P_r}} )=
{\cal C} _* ({\overline {P_r}}, Q_r )$.

\paragraph{Now we can finally start the proof.}
$\overline {P_r}$ is a compact $C^1$ manifold with boundary, and 
$A'_r,\,\,B'_r$ are its
$C^1$ submanifolds with boundaries, which are relative cycles in 
${\cal C} _* ({\overline {P_r}}, \partial {\overline{P_r}} )$.

Now we have to choose $\varepsilon _1$ and $K$ properly. Observe that if 
a cube, having a face with vertices $v_1,v_2,v_3,v_4$, in this cyclic order, is
inscribed in $f(S^2)$, then
$$\min_{u \in S^2 } \|f(u)\| \le \|v_1\|=...=\|v_4\| \le \max_{u \in S^2}
\|f(u)\|,$$
moreover, for $1 \le i < j \le 4$, we have
$\|v_i \pm v_j\| \in \{ 2\|v_i\|/
\sqrt 3, 2\sqrt 2 \|v_i\|/\sqrt 3 \}$, so $\|v_i \pm v_j\| /\sqrt {2} \ge 
\|v_i\| \sqrt {2/3}$. 
We have 
to choose $\varepsilon _1 $ so small and $K$ so large that, taking 
${\overline {P_r}} \subset X/G$
rather than $X/G$ does not exclude any inscribed cube. We have 
$${\overline {P_r}} \supset P_r 
\supset (X_ r /G) \setminus \left(  \varphi ^{-1} \left( (1, 
5)\right) /G \right) \supset$$
$$\left(  \varphi ^{-1} \left( [0, 
1] \right) /G \right) \setminus \left(  \varphi ^{-1} \left( (1, 
5) \right) /G \right) 
=  \varphi ^{-1} \left( [0, 
1] \right) /G .$$
Hence for $x \in X$ and $x/G \not\in \overline{P_r}$ we have $\varphi (x) >1$, 
thus either $\|v_i\| < \varepsilon _1 $ for some $1 \le i \le 4$, or $\|v_i\|
>K$ for some $1 \le i \le 4$, or $\|v_i \pm v_j\|/\sqrt {2} < \varepsilon _2$ 
for some $ 1 \le i < j \le 4$. However, for $v_1,...,v_4$ being vertices of a 
cube, like above, we have in the last case $\|v_i\| \le \left( \|v_i \pm v_j\| 
/ \sqrt {2} \right) \sqrt {3/2} < \varepsilon _2 \sqrt {3/2} < \varepsilon _1$,
by the choice of $\varepsilon _2 $. Hence the first or second inequality must 
hold. These however can be excluded by choosing 
$$0 < \varepsilon _1 < \min_{u \in S^2} \|f(u)\| \le \max_{u \in S^2} \|f(u)\|
<K.$$
Thus in fact taking $\overline {P_r} \subset X/G $ rather than $X/G$ does not 
exclude any inscribed cube.
 

Now recall that $A'_r$ and $B'_r$ are $C^1$ manifolds with boundaries, 
with their boundaries lying in $ \partial {\overline{P_r}} =Q_r$. Let 
$$Q_{r1}:=Q_r \cap \left[ \left\{ (v_1,...,v_4) \in X \mid \max \{\frac 
{\varepsilon _1}{\|v_i\|}, \frac {\|v_i\|}{K} \mid 1 \le i \le 4 \} <
\right. \right. $$
$$\left. \left. \max \{\frac {\sqrt {2} \varepsilon _2}{\|v_i \pm v_j\|} \mid 
1 \le i < j \le 4 \} \right\} \biggm/ G \right]$$
and 
$$Q_{r2}:=Q_r \cap \left[ \left\{ (v_1,...,v_4) \in X \mid \max \{\frac 
{\varepsilon _1}{\|v_i\|}, \frac {\|v_i\|}{K} \mid 1 \le i \le 4 \} > 
\right. \right. $$
$$\left. \left. \max \{\frac {\sqrt {2} \varepsilon _2}{\|v_i \pm v_j\|} \mid 
1 \le i < j \le 4 \} \right\} \biggm/ G \right].$$
Then $Q_{r1}$ and $Q_{r2}$ are disjoint open subsets of $Q_r$.
Actually we will attain that $A'_r$ will be a 
relative cycle in ${\cal C} _*(\overline {P_r}, Q_{r2})$, 
and $B'_r$ will be a 
relative cycle in ${\cal C} _*(\overline {P_r}, Q_{r1})$, 
i.e., $\partial A'_r \subset
Q_{r2}$ 
and  $\partial B'_r \subset
Q_{r1}$. That is, we 
want to attain $(\partial A'_r) \cap (Q_r \setminus Q_{r2})=\emptyset $ and 
$(\partial B'_r) \cap (Q_r \setminus Q_{r1})=\emptyset $. 
However, by $\partial A'_r, 
\partial B'_r \subset Q_r \subset \varphi ^{-1} \left( [1,5] \right)/G$ we have
$$\partial A'_r \subset \left[ A/G \right] \cap \left[ \varphi ^{-1} \left(
[1,5] \right) /G \right] {\rm {and}} \,\,
\partial B'_r \subset \left[ f(S^2)^4/G \right] \cap \left[ \varphi ^{-1} 
\left( [1,5] \right) /G \right].$$
Therefore, it suffices to establish
$$ \left[ A/G \right] \cap \left[ \varphi ^{-1} \left(
[1,5] \right) /G \right] \cap (Q_r \setminus Q_{r2}) 
= \emptyset \,\,{\rm {and}} 
\left[ f(S^2)^4/G \right] \cap \left[ \varphi ^{-1} 
\left( [1,5] \right) /G \right] \cap (Q_r \setminus Q_{r1})=\emptyset.$$
To obtain a contradiction, we assume that one of these sets is nonempty.

First assume that $(v_1,...,v_4)/G \in \left[ A/G \right] \cap \left[ \varphi
^{-1} \left( [1,5] \right) /G \right] \cap (Q_r \setminus Q_{r2})$. 
Then $1 \le \varphi 
(v_1,...,v_4) \le 5 $, hence, by taking in consideration the definitions of 
$\varphi $ and $Q_{r2}$, we have 
$$\max \left\{ \frac 
{\varepsilon _1}{\|v_i\|}, \frac {\|v_i\|}{K} \mid 1 \le i \le 4 \right\} 
\in [0,5]
\,\,{\rm {and}}\, 
\max \left\{ \frac {\sqrt {2} \varepsilon _2}{\|v_i \pm v_j\|} \mid 1 
\le i < j \le 4 \right\} \in [1,5].$$
These imply $\|v_i\| \ge \varepsilon _1 /5 $ for each $1 \le i \le 4$, and 
$\|v_i \pm v_j\| / \sqrt {2} \le \varepsilon _2 $ for some $1 \le i < j \le 
4 $ and some choice of the sign $\pm $. However, for vertices of a cube 
$v_1,...,v_4 $ like above, we have $\|v_i \pm v_j \| / \sqrt {2} \ge 
\|v_i \| \sqrt {2/3}$, that implies for the above pair $\{i,j\}$ and the above 
choice of the sign that 
$$\frac {\varepsilon _1 }{100} = \varepsilon _2 \ge \frac {\| v_i \pm v_j \|}
{\sqrt {2}} \ge \|v_i \| \sqrt {\frac {2}{3}} \ge \frac {\varepsilon _1 }{5}
\sqrt {\frac {2}{3}},$$
a contradiction.

Second assume that $(v_1,...,v_4)/G \in \left[ f(S^2)^4/G \right] \cap 
\left[ \varphi ^{-1} \left( [1,5] \right) /G \right] \cap (Q_r \setminus 
Q_{r1}) $. Then $1 
\le \varphi (v_1,...,v_4) \le 5$, hence, by the definitions of $\varphi $ 
and $Q_{r1}$, we have 
$$\max \left\{ \frac 
{\varepsilon _1}{\|v_i\|}, \frac {\|v_i\|}{K} \mid 1 \le i \le 4 \right\} 
\in [1,5].
$$
Then for some $1 \le i \le 4 $ we have $\|v_i\| \le \varepsilon _1 $ or 
$\|v_i\| \ge K$. However, by $(v_1,...,v_4)/G \in f(S^2)^4 /G $ any of these 
possibilities contradicts the choice of $\varepsilon _1 $ and $K$. This ends
the proof of the statement that $A'_r$ or $B'_r$ are relative cycles in 
${\cal C} _*(\overline{P_r},Q_{r2})$ or ${\cal C} _*(\overline{P_r},Q_{r1})$,
respectively.

Actually, we shall choose the numbers $\varepsilon _1 $ and $K$ 
even more carefully
as follows. We define $g$ to be the embedding of the ellipsoid 
as in the plan of the proof.
By genericity of $Q_r$ we may assume that $Q_r$ intersects 
the orbifold $g(S^2)^4/G$ (that is a manifold in a neighbourhood of $Y_r /G$)
transversally --- and hence in finitely many points --- as well. 
By hypothesis we have an odd homotopy $H$ (i.e., one satisfying $H(-u,t)=
-H(u,t)$ for $u \in S^2,\,\,t \in [0,1]$) 
between $f$ and the standard embedding of $S^2$ to $\R^3$. Similarly 
there is an odd homotopy, obtained by linear interpolation, between the 
standard embedding of $S^2$ to $\R^3$ and $g$. Putting these together we 
obtain an odd homotopy $H_1:S^2 \times [0,1] \to {\R}^3 \setminus \{0\}$ 
between $f$ and $g$. Now we suppose 
$$0 < \varepsilon _1 < \min \{\|H_1 (u,t)\| \mid u \in S^2,\,\,t \in [0,1] \} 
\le \max \{\|H_1 (u,t)\| \mid u \in S^2,\,\,t \in [0,1] \} <K.$$

Now $H_1^4:(S^2)^4 \times [0,1] \to ({\R}^3 \setminus \{0\} )^4,\,\,
H_1^4((u_1,...,u_4),t)=(H_1(u_1,t),...,H_1(u_4,t))$ is a homotopy (with four odd 
components) betwen $f \times f \times f \times f = f^4:(S^2)^4 \to ({\R}^3 
\setminus \{0\})^4$ and $g^4:(S^2)^4 \to ({\R}^3 
\setminus \{0\})^4$. Moreover, $H_1^4$ is $G$--equivariant, i.e., for $T\in G 
$ we have $H_1^4(T(u_1,...u_4),t)=T H_1^4 ((u_1,...,u_4),t)$ (this being true 
even for $T \in S_4 \times {\Z}_2^4$, the $i$-th ${\Z}_2$ acting on the $i$-th 
coordinate by $v \mapsto \pm v$, and on the other coordinates identically). 
Therefore $H_1^4$ induces a homotopy $H_1^4/G:(S^2)^4 \times [0,1])/G \to
({\R}^3 \setminus \{0\})^4/G$ ($G$ acting on the factor $[0,1]$ trivially) 
between $f^4/G,g^4/G :(S^2)^4 /G \to
({\R}^3 \setminus \{0\})^4/G$, by the formula $(H_1^4/G)(G(u_1,...,u_4),t)=
GH_1^4((u_1,...u_4),t)$, where we put the quotient topology on each of these
spaces. (There is one point to be clarified: the restriction of the quotient 
topology ${\cal T}=((S^2)^4 \times [0,1])/G$ to $((S^2)^4 \times \{0\})/G$ 
is the quotient topology ${\cal T}_0=(S^2)^4/G$,
and a similar statement with $\{1\}$ rather than $\{0\}$. 
This follows from ${\cal T}
={\cal T}_0 \times [0,1]$. To see this last statement, observe that we have 
an evident mapping $(S^2)^4 \times [0,1] \to ((S^2)^4 / G) \times [0,1]$. This 
factors through the quotient topology ${\cal T}$. However ${\cal T}$ is 
compact, and $((S^2)^4 / G) \times [0,1]$ is $T_2$, so a bijective map 
between them is a homeomorphism. The $T_2$ property of $((S^2)^4 / G) \times 
[0,1]$ 
follows since the quotient space of a $T_2$ space by a finite group of 
homeomorphisms is a $T_2$ space.)

Now recall that $Q_r$ is transversal to $f(S^2)^4/G$ and $g(S^2)^4/G$, which 
are manifolds in a neighbourhood of $Q_r$. In other words, the maps $f^4/G$ 
and $g^4/G$ are transversal to $Q_r$. Therefore also the homotopy $H_1^4/G$
connecting them is transversal to $Q_r$ at $t=0$ and $t=1$.

Our aim will be to modify the map $H_1^4/G:((S^2)^4 / G) \times [0,1] \to 
({\R}^3 \setminus \{0\})^4 / G $ a bit, so as to obtain a map transversal 
to $Q_r$. For this aim first observe that $(S^2)^4/G$ is not a manifold. 
Therefore, analogously to $\varphi , \psi :X=({\R}^3)^4 \setminus S \to {\R}$ 
we define $\Phi , \Psi :(S^2)^4 \setminus S \to {\R}$, by the formulas 
$$\Phi (u_1,...u_4) := \max \left\{ \frac {\sqrt 2 \delta}{\|u_i \pm u_j \|}
\mid 1 \le i < j \le 4 \right\}$$ 
and
$$\Psi (u_1,...,u_4) :=\left( \sum _{1 \le i < j \le 4} \frac {2 \delta ^2}
{\|u_i \pm u_j\|^2} \right)^{1/2},$$
where $\delta >0$ is to be chosen later to be sufficiently small. We have 
$$\Phi \le \Psi \le {\sqrt {12}}\Phi .$$
Let $s \in ({\sqrt {12}},4)$ be a not critical value of the $C^{\infty}$
function $\Psi $, and let 
$$M_s :=\Psi ^{-1}([0,s))=\Psi ^{-1} ((0,s)) \supset \Phi ^{-1}([0,1]) 
{\rm {\,\,and\,\,}} N_s := \Psi ^{-1} (s) \subset \Phi ^{-1} ((1,4)).$$
Then $M_s \subset \Phi ^{-1} ([0,4))$ and $N_s = \partial {\overline {M_s}}$ 
is a 
non--empty closed embedded $C^{\infty}$ $7$--submanifold of $(S^2)^4 
\setminus S$. Moreover $M_s$ and $N_s$ are $G$--invariant. The factor space 
${\overline {M_s}}/G=(M_s \cup N_s) /G$ is a compact $C^{\infty }$ manifold 
with boundary $N_s / G $.

Observe that the odd homotopy $H_1:S^2 \times [0,1] \to {\R}^3 \setminus 
\{0\}$ has a compact domain, hence is uniformly continuous. Moreover we have 
$H_1^4\left( ((S^2)^4 \setminus M_s) \times [0,1] \right) \subset 
H_1^4\left( ((S^2)^4 \setminus \Phi ^{-1}([0,1])) \times [0,1] \right) $.
Now choose $\delta $ so that $u_i,u_j \in S^2,\,\,\|u_i-u_j\|/\sqrt {2} 
< \delta $ implies 
$\| H_1(u_i,t)-H_1(u_j,t)\| /{\sqrt 2} < \varepsilon _2 /7$. 
By oddness then we 
have that also $u_i,u_j \in S^2,\,\,\|u_i+u_j\|/\sqrt {2} < \delta $ implies 
$\| H_1(u_i,t)+H_1(u_j,t)\| /{\sqrt 2} < \varepsilon _2 /7$. Then we claim
$$H_1^4\left( ((S^2)^4 \setminus M_s) \times [0,1] \right) \subset 
H_1^4\left( ((S^2)^4 \setminus \Phi ^{-1}([0,1])) \times [0,1] \right)
\subset \varphi ^{-1} ((7, \infty )) \cup S.$$
In fact, for this it suffices to show that $u_1,...,u_4 \in S^2,\,\,\Phi (u_1,
...,u_4)>1$ and $(H_1(u_1,t),...,$
\newline $H_1(u_4,t)) \not\in S$ imply $\varphi 
(H_1(u_1,t),...,H_1(u_4,t))>7$. We have $\Phi (u_1,...,u_4) >1 
\Longleftrightarrow \min _{i<j} \|u_i \pm u_j \| /{\sqrt 2} < \delta $, hence 
at these conditions $\varphi (H_1(u_1,t),...,H_1(u_4,t)) \ge \max _{i<j}$
\newline $ ({\sqrt 2}
\varepsilon _2 )/\| H_1(u_i,t) \pm H_1(u_j,t)\| >7$, 
by the choice of $\delta $ and oddness of the homotopy $H_1$. 

We have 
$Q_r \subset \varphi ^{-1}([1,5])/G \subset V_1:= \varphi ^{-1}((0,6))/G 
\subset {\overline {V_1}} \subset V_2:= \varphi ^{-1}((0,7))/G$. By $H_1^4
\left( ((S^2)^4 \setminus M_s ) \times [0,1] \right) \subset 
\varphi ^{-1} ((7, 
\infty )) \cup S$ we have $\left( ((S^2)^4 \setminus M_s ) \times [0,1] \right)
\cap (H_1^4)^{-1} $
\newline $\left( \varphi ^{-1}((0,7)) \right) = \emptyset $, i.e., 
$(H_1^4)^{-1} \left( \varphi ^{-1}((0,7)) \right) \subset M_s \times [0,1]$. Here 
$M_s$ is a $C^1$ manifold, thus $M_s \times [0,1]$ is a $C^1$ manifold with 
boundary. We have $(H_1^4/G)^{-1} \left( \varphi ^{-1}((0,7))/G \right)
=(H_1^4/G)^{-1}(V_2) \subset (M_s/G) \times [0,1]$, so $(H_1^4/G)^{-1}(V_2)$
is an open set in the $C^1$ manifold with boundary $(M_s/G) \times [0,1]$, so 
it is also a $C^1$ manifold with boundary. The set $(H_1^4)^{-1}(V_2 \setminus 
V_1)$ is a relatively closed set in $(H_1^4/G)^{-1}(V_2)$, and on it the map 
$H_1^4/G$ is transversal to $Q_r$, since $(H_1^4/G)[(H_1^4/G)^{-1}(V_2 
\setminus V_1)] \subset 
V_2 \setminus V_1$ is disjoint to $Q_r$. Moreover, as we have 
seen above, also on $[((S^2)^4/G) \times
\{0,1\}] \cap (H_1^4/G)^{-1}(V_2)$ the 
map $H_1^4/G$ is transversal to $Q_r$. 

Hence by \cite{brocker-janich}, 
Theorem 14.7 there exist arbitrarily small
$C^1$--perturbations $H'$ of $H_1^4/G$ on the set $(H_1^4/G)^{-1}(V_2)$ such 
that $H'$ is transversal to $Q_r$, and $H'$ coincides with $H_1^4/G$ on a 
neighbourhood of the closed set $(H_1^4/G)^{-1}(V_2 \setminus V_1) \cup
\left\{ [((S^2)^4/G) \times \{0,1\}] \cap (H_1^4/G)^{-1} \right.$
\newline $\left. (V_2) \right\} $.
Since the compact range of $H_1$ lies in the set $\{v \in {\R}^3 \setminus 
\{0\} \mid \varepsilon _1 < \| v \| <K\}$, therefore we may assume that $H'
[(H_1^4/G)^{-1}(V_2)] \subset \{ (v_1,...,v_4) \in ({\R}^3 \setminus \{0\})^4
\mid \varepsilon _1 < \| v_i \| <K\}/G$. We can define $H'$ also on $(H_1^4/G)
^{-1}[((\R^3 \setminus \{0\})^4/G)
\setminus (\overline {V_1})]$, as $H_1^4/G$. Then the last inclusion formula about the range of 
$H'$ remains valid, in particular $H'$ has a range in $({\R}^3 \setminus 
\{0\})^4/G$. Thus we have defined $H'$ on two open sets, coherently on their 
intersection, which will together give a perturbation $H'$ of $H_1^4/G$, 
transversal to $Q_r$. In particular, $H'^{-1}(Q_r)$ is a manifold (here only 
$H'| \left( (H_1^4/G)^{-1}(V_2) \right)$ plays a role). Moreover, also 
$\left( (S^2)^4/G \times \{0,1\} \right) \cap H'^{-1}(Q_r)$ is a manifold.
Namely, $Q_r$ is transversal to $f(S^2)^4/G$ and $g(S^2)^4/G$, i.e., $f^4/G,
g^4/G :(S^2)^4/G \to ({\R}^3 \setminus \{0\})^4/G$ are transversal to $Q_r$, 
so $(f^4/G)^{-1}(Q_r),(g^4/G)^{-1}(Q_r) \subset (S^2)^4/G$ are manifolds. 
Further, $H'|\left((S^2)^4/G \times \{0\} \right) =(H_1^4/G)|
\left((S^2)^4/G \times \{0\} \right)=f^4/G$ and 
$H'|\left((S^2)^4/G \times \{1\} \right) = g^4/G$.

Let us consider the restriction $H'|H'^{-1}({\overline {P_r}})$. We 
assert that we may assume that $H'^{-1}({\overline {P_r}})$ is a subset of 
the manifold with boundary $(M_s/G) \times [0,1]$, or in other words, $H'^{-1}
\left( (({\R}^3)^4/G) \setminus {\overline {P_r}} \right) \supset \left( ((S^2
)^4 \setminus M_s)/G \right) \times [0,1]$, or in yet other words, $H'\left\{ 
\left( ((S^2)^4 \setminus M_s)\right. \right. $
\newline $\left. \left. 
/G \right) \times [0,1] \right\} \subset (({\R}^3)^4 /G) 
\setminus {\overline {P_r}}$. In fact, we have from above $(H_1^4/G)\left\{
\left( ((S^2)^4 \setminus M_s)/G \right) \right. $
\newline $\left. \times [0,1] \right\} \subset (\varphi ^{-1}((7,
\infty )) \cup S)/G) $ 
and also ${\overline {P_r}}=P_r \cup Q_r \subset \varphi 
^{-1}([0,5])/G$, hence $\left( \varphi ^{-1}((5, \right. $
\newline $\left. \infty )) \cup S \right) /G 
\subset (({\R}^3)^4/G) \setminus {\overline {P_r}}$.
Now observe that $\left( \varphi ^{-1}((5, \infty )) \cup S \right) /G$ is a 
neighbourhood of the closed set 
$\left( \varphi ^{-1}([7, \infty )) \cup S \right) /G$ containing the open set
$\left( \varphi ^{-1}((7, \infty )) \cup S \right) /G$. The complement of this 
neighbourhood is a compact set, hence this neighbourhood contains even some
metric $\eta $--neighbourhood of 
$\left( \varphi ^{-1}([7, \infty )) \cup S \right) /G$. Thus a sufficiently small 
$C^1$--perturbation $H'$ of $H_1^4/G$  satisfies 
$$H'\left\{ \left( ((S^2)^4 \setminus M_s) /G \right) \times [0,1] \right\} 
\subset \left( \varphi ^{-1}((5, \infty )) \cup S \right) /G \subset (({\R}^3)^4
/G) \setminus {\overline {P_r}}.$$

Moreover, $H'^{-1}({\overline {P_r}})\,\,( \subset (M_s /G) \times [0,1])$ 
is itself 
a $C^1$--manifold with boundary. Namely, we consider $H'|\left( 
(M_s /G) \times [0,1] \right):(M_s /G) \times [0,1] \to ({\R}^3 \setminus 
\{0\})^4/G$. Here $(M_s /G) \times [0,1]$ is a manifold with boundary 
$(M_s /G) \times \{0,1\}$, $({\R}^3 \setminus 
\{0\})^4/G$ is an orbifold, containing the manifold with boundary $
{\overline {P_r}} \subset \varphi ^{-1}([0,5])/G$, that is contained in the 
manifold $\varphi ^{-1}([0,6))/G=\varphi ^{-1}((0,6))/G$. Then $[H'|\left( 
(M_s /G) \times [0,1] \right)]^{-1} \left( \varphi ^{-1}((0,6))/G \right)$ 
is an open subset of
$(M_s /G) \times [0,1]$, hence is a manifold with boundary, and contains
$[H'|\left( 
(M_s /G) \times [0,1] \right)]^{-1} ({\overline {P_r}})=H'^{-1}
({\overline {P_r}})$. Now the restriction $H'' $ of $H'$ to $[H'|\left( 
(M_s /G) \times \right. $
\newline $\left. [0,1] \right)]^{-1} \left( \varphi ^{-1}((0,6))/G \right)$
maps this manifold with boundary to the manifold $\varphi ^{-1}((0,6))/G$.
Then by [Hir], Ch. 1, Theorem 4.2, (ii) the inverse image by this map of the 
manifold with boundary ${\overline {P_r}}$, i.e., $H'^{-1}({\overline {P_r}})$,
is itself a manifold with boundary. In fact, for this we need that 
\newline (1) $H''$ is 
transversal to $Q_r$, and 
\newline (2) its restriction to the boundary of its domain, 
i.e., to $[H'|\left( 
(M_s /G) \times \{0,1\} \right)]^{-1} \left( \varphi ^{-1}((0,\right. $
$\left. 6))/G \right)$
is also transversal to $Q_r$. 
\newline Here (1) follows, since 
$H''$ pointwise coincides 
with $H'$, and $H'$ is transversal to $Q_r$. Moreover (2) follows, since on 
$(M_s /G) \times \{0,1\}$ (even on $(S^2)^4/G \times \{0,1\}$) $H''$ and $H'$ 
coincide with $H_1^4/G$, and from above $H_1^4/G$ is transversal to $Q_r$ at 
$t=0$ and $t=1$. 

We make each calculation (boundary, homology, etc.) mod 2. Then we have 
$$\partial (H'|H'^{-1}({\overline {P_r}}))=H'|[((S^2)^4/G \times \{0,1\}) \cap 
H'^{-1}({\overline {P_r}})]+H'|H'^{-1}(Q_r),$$
where $((S^2)^4/G \times \{0,1\}) \cap 
H'^{-1}({\overline {P_r}}) \subset ((S^2)^4/G \times \{0,1\}) \cap ((M_s/G)
\times [0,1]) = (M_s/G) \times \{0,1\}$, and $H'^{-1}(Q_r)$ is a manifold, 
from above. Further, $((S^2)^4/G \times \{0,1\}) \cap 
H'^{-1}({\overline {P_r}})$ is a manifold with boundary. This can be seen 
analogously like above. We will show this for $((S^2)^4/G \times \{0\}) \cap 
H'^{-1}({\overline {P_r}})$. (The other case is analogous.) Then $ 
{\overline {P_r}}$ is a subset of the manifold 
$\varphi ^{-1}((0,6))/G$, and we 
have to investigate the set $(f^4/G)^{-1}({\overline {P_r}}) \,\,( \subset
M_s/G)$. The restriction of $f^4/G$ to $M_s/G$ maps the manifold $M_s/G$ to 
$({\R}^3 \setminus \{0\})^4/G$. We consider the inverse image of 
$\varphi ^ {-1} 
((0,6))/G$ by this map, that is a manifold contained in $M_s/G$, and we 
further restrict the above restriction of $f^4/G$ to it, obtaining a map to 
$\varphi ^{-1} ((0,6))/G$. By \cite{hirsch}, Ch. 1, Theorem 4.2, (ii) 
we have that $(f^4/G)^
{-1}({\overline {P_r}})$ is a manifold with boundary, provided that our map 
is transversal to $Q_r$. However, our map pointwise coincides with $f^4/G$, so 
this follows from the choice of the manifold $Q_r$.

By the above paragraph the cycle 
$$H'|[((S^2)^4/G \times \{0,1\}) \cap H'^{-1}(
{\overline {P_r}})]+H'|H'^{-1}(Q_r) =$$
$$(f^4/G)|[
(f^4/G)^{-1}({\overline {P_r}})]+(g^4/G)|[
(g^4/G)^{-1}({\overline {P_r}})]+$$
$$H'|H'^{-1}(Q_r)$$ 
is homologous to $0$.

Now we assert that the range of the map $H'|H'^{-1}(Q_r)$ is included in 
$Q_{r1}$. It will suffice to prove $H'\left( (S^2)^4/G \times [0,1] \right) 
\cap (Q_r \setminus Q_{r1}) \subset \left( \{ (v_1,...,v_4) \mid \varepsilon _1
 < \| v_i \| <K \}/G \right) \cap (Q_r \setminus Q_{r1}) = \emptyset $. Suppose
 that $(v_1,...,v_4)/G$ belongs to the last set. Then $\max _i \{ 
 \varepsilon _1 /\| v_i \|, \|v_i \| /K \} $
\newline $ <1$, and, by the definition of 
 $Q_{r1}$, also $\max _{i<j}\{({\sqrt 2} \varepsilon _2)/ \| v_i \pm v_j \| \} 
 <1$, hence $\varphi (v_1,...,v_4) <1$. However, $Q_r \setminus Q_{r1} \subset 
 Q_r \subset \varphi ^{-1} ([1,5])/G$, so $\varphi (v_1,...,v_4) \ge 1$, a 
 contradiction.

Then by the above homological formula $(f^4/G)|[
(f^4/G)^{-1}({\overline {P_r}})]$ and $(g^4/G)|[
(g^4/G)^{-1}({\overline {P_r}})]$ are homologous, mod $Q_{r1}$. Now recall 
that $A_r'=[A/G] \cap {\overline {P_r}}$ is a relative cycle in 
$({\overline {P_r}}, Q_{r2})$, and 
$B_r'=[f(S^2)^4/G] \cap {\overline {P_r}}$ is a relative cycle in 
$({\overline {P_r}}, Q_{r1})$. 
Since the roles of $f$ and $g$ are symmetric, also 
$[g(S^2)^4/G] \cap {\overline {P_r}}$ is a relative cycle in 
$({\overline {P_r}}, Q_{r1})$. By the
above homologicity property the ranges of the injections $(f^4/G)|
[
(f^4/G)^{-1}({\overline {P_r}})]$ and $(g^4/G)|[
(g^4/G)^{-1}({\overline {P_r}})]$, i.e., $[f(S^2)^4/G] \cap 
{\overline {P_r}}$ and $[g(S^2)^4/G] \cap 
{\overline {P_r}}$ will have the same intersection numbers (mod $2$) with the 
relative cycle $A_r'$.

Therefore the intersection numbers of the relative cycle $\alpha $ 
realized by the set $A'_r$ -- which is a
relative cycle in $({\overline {P_r}},
Q_{r2})$ --  with $\beta _f$, realized by the set $B_r'$,
and $\beta _g$, realized analogously by the set $[g(S^2)^4/G] \cap
{\overline {P_r}}$ -- which in turn are relative cycles in $({\overline {P_r}}
,Q_{r1})$ -- will be equal. 
$\square$

\medskip
(Alternatively one can give the same proof in the language of cohomologies 
by turning to the dual
cohomology classes. 
Let  $[\alpha ]$ be the homology class of $\alpha $ in $H_*({\overline {
P_r}}, Q_{r2};
\Z_2)$, $[\beta _f] = [\beta _g]$ the homology class of 
$\beta _f$ 
(and of $\beta _g$) in $H_*({\overline {
P_r}}, Q_{r1};\Z_2)$.
Their Poincar\'e duals will be denoted by 
$$D[\alpha ] \in H^*({\overline {P_r}}, 
Q_{r1};\Z_2)\
\ {\rm and}\
\ D[\beta _g] = D[\beta _f] \in H^*({\overline {P_r}}, Q_{r2};\Z_2).$$
Their product  $D[\alpha ] \cup D[\beta _g] $ can be considered as an element 
in $H^{{\rm {dim}}{\overline {P_r}}}(
{\overline {P_r}}, \partial {\overline {P_r}}; \Z_2)
= \Z_2$. Thus we have 
$$\ {\rm the\,\,number\,\,of\,\,cubes\,\,on\,\,} f(S^2) = D[\alpha 
\cap \beta _f] = D[\alpha ] \cup D [\beta _f]=$$
$$ D[\alpha ] \cup D[\beta _g] = D[\alpha  \cap \beta _g] 
= \ {\rm the\,\,number\,\,of\,\,cubes\,\,on\,\,} g(S^2)
\ne 0.)$$
$\square$

Actually the same arguments yield the following more general theorem (compare 
Theorem C of \cite{griffiths}, cited at the beginning of this section).

\begin{theorem}
\label{konkav'} 
Let $f$ be as in Theorem~\ref{konkav}. Then there is a square
based box in ${\R}^3$, with its centre in the origin, and with any given
ratio of the height to the basic edge, having all its vertices on the surface
$f(S^2)$.
\end{theorem}

\begin{proof}
We proceed identically as above, with the only difference that we apply
Corollary~\ref{ellipsoid} not for a cube, but use its consequence
mentioned in the proof of Proposition~\ref{d4}, again observing that $3$ is
odd.
\end{proof}

\paragraph{Question.} Can one generalize further Theorem~\ref{konkav'}, for any
box, with any given ratios of edge-lengths (like in Theorem~\ref{boxes})?
(Observe, that here Corollary~\ref{ellipsoid} is of no use, since 
 $3!/(1!)^3=6$ is even. Our proof fails for this case, because the cycles,
whose intersection points give the cubes on the surface, were not oriented,
i.e., integer cycles, because of the presence of orientation reversing
elements in the group $G$.) 

\paragraph{Remark.} Our Theorem~\ref{konkav} generalizes Theorem~\ref{cube} and
Theorem~\ref{konkav'} generalizes Theorem~\ref{square}. In fact, for $ 
\min F>0$ we may choose the function $f(u)$ from Theorems~\ref{konkav} and
~\ref{konkav'} as $F(u) \cdot u$, where $F(u)$ is the function from
Theorems~\ref{cube} and ~\ref{square}. This shows Theorems~\ref{cube} and
~\ref{square} for even $C^1$ functions $F$ (the condition $\min F >0$ can
evidently be suppressed). To obtain their statements for
even continuous functions $F$, we use density of $C^1(S^2)$ in $C^0(S^2)$, and
compactness of $SO(3)$.


\section{Universal covers: proof of Theorem~\ref{rhombic}}
\label{covers}

As we have explained in Section~\ref{introduction}, P\'al in \cite{pal}
proved that a regular hexagon with distance $1$ between its opposite
sides was a universal cover in $\R^2$ . In this section 
we deal with Makeev's generalization
of P\'al's question. 

Let $\Sigma ^n \subset {\R}^n$ be a regular simplex of edge-length $1$,
with vertices $v_1,...,v_{n+1}$. Let $U_n$ be the intersection of $n(n+1)/2$ 
parallel strips $S_{ij}\,\,(1 \le i <j \le n+1)$ of width 1, where $S_{ij}$ is 
bounded by the 
$(n-1)$-planes orthogonal to the segment $[v_i,v_j]$, and passing through 
$v_i$ and $v_j$, respectively. For $n=2$ we have that $U_2$ is the above 
mentioned regular hexagon. Makeev in \cite{makeev3} proposes, as a
generalization of P\'al's result, the
following

\begin{conjecture}[Makeev] $U_n$ is a universal cover 
in ${\R}^n$.
\end{conjecture} 

First we reformulate this conjecture, following \cite{makeev0}, 
to a problem concerning
continuous functions.

Let $K \subset {\R}^n$ be a non-empty compact convex set, of diameter at 
most $1$. One can define its support function $h:S^{n-1} \to {\R}$ 
as in \cite{bonnesen-fenchel}. 
Let $A \in SO(n)$ and $1 \le i 
< j \le n+1$, and consider the parallel 
supporting $(n-1)$-planes 
$$\{ x \in
{\R}^n \mid \langle x, A(v_j-v_i) \rangle = h(A(v_j-v_i)) \}$$ and $$ \{ x \in
{\R}^n \mid \langle x, A(v_i-v_j) \rangle = h(A(v_i-v_j)) \}$$ of $K$, 
whose distance is at most $1$, because $K$ has diameter at most $1$. 
Consider also their mid-$(n-1)$-plane $$\{x \in {\R}^n
\mid \langle x, A(v_j-v_i) \rangle = [h(A(v_j-v_i))-h(A(v_i-v_j))]/2\}.$$ 
Let these supporting $(n-1)$-planes bound the parallel strip $S_{ij}'$. 
Suppose that
the above mid-$(n-1)$-planes are concurrent, with (unique) common point $x 
\in {\R}^n$. Then include each strip $S_{ij}'$ to a strip of width $1$,
having the same mid-$(n-1)$-plane as $S_{ij}'$. These larger strips are of 
the form $AS_{ij}+x$. That is, $K \subset \cap S_{ij}' \subset A(\cap S_{ij})
+x=AU_n+x$, and thus $U_n$ would be a universal cover. Therefore 
Makeev's conjecture would follow from 

\begin{conjecture}
Let $F:S^{n-1} \to {\R}$ be an odd function, and
let $\Sigma _n \subset {\R}^n$ be a regular simplex of edge-length $1$, 
with vertices $v_1,...,v_{n+1}$. Then there exists an $A \in SO(n)$ such that 
the $n(n+1)/2\quad
(n-1)$-planes $$\{x \in {\R}^n \mid \langle x, A(v_j-v_i) \rangle =
F(A(v_j-v_i))\}\,\,(1 \le i <j \le n+1)$$ are concurrent. 
\end{conjecture}

We remark that \cite{makeev0} 
has given this conjecture in a slightly different 
formulation.

Observe that $F$ induces a map $\Phi :SO(n) \to {\R}^{n(n+1)/2}$, of 
coordinates $F(A(v_j-v_i)),\,\,1 \le i < j \le n+1$. We need $A \in SO(n)$, 
such that $\Phi (A)$ lies in the $n$-subspace 
\begin{eqnarray*}
\{( \langle x, A(v_j-v_i) 
\rangle )_{1 \le i < j \le n+1} \mid x \in {\R}^n \}& =& \{ (\langle A^*x,
v_j-v_i \rangle )_{1 \le i < j \le n+1} \mid x \in {\R}^n \}\\ 
&=&\{ (y, v_j-v_i
\rangle )_{1 \le i < j \le n+1} \mid y \in {\R}^n \}
\end{eqnarray*} of ${\R}^{n(n+1)
/2}$. Note that ${\rm dim} SO(n)=n(n+1)/2-n$, hence the number of variables
equals the number of constraints, and thus there is a chance to apply 
intersection theory. Analogously as in \S 2, it would be sufficient to find a 
test odd function $F_0:S^{n-1} \to {\R}$, for which we have a unique 
solution $A \in SO(n)$ (up to symmetries of the convex 
polyhedron $U_n$), and for
which the intersection number at this unique solution is non-zero. A 
candidate for $F_0$ could be e.g. $(h_0(u)-h_0(-u))/2$, where $h_0$ is the 
support function of $\Sigma ^n$, for which $A=I$ is a solution. (For $A=I$
the mid-$(n-1)$-planes are just the halving $(n-1)$-planes of the edges of 
$\Sigma ^n$, all passing through the centre of $\Sigma ^n$.) For $n=2$ 
elementary considerations show unicity of the solution (up to symmetries of 
the regular hexagon), and in a natural sense the intersection is transversal 
there.  

Though we could not complete the above program, for $n=3$ we can 
deduce Makeev's conjecture from Proposition~\ref{s4}:  

\begin{theorem} Let $F:S^2\rightarrow \R$ be an odd function, and let 
$\Sigma_3\subset \R^3$ be a regular simplex of edge-length 1 centred at
the origin, with vertices 
$v_1,v_2,v_3$ and $v_4$. Then there exists an $A\in SO(3)$ such that the 
$6$ planes  $$\{x \in {\R}^3 \mid \langle x, A(v_j-v_i) \rangle =
F(A(v_j-v_i))\}\,\,(1 \le i <j \le 4)$$ are concurrent.
\label{concurrent}
\end{theorem}

\begin{proof} 
Similarly as above  $F$ induces a map $\Phi :SO(3) \to \R^6$ given
by  coordinates $F(A(v_j-v_i)),\,\,1 \le i < j \le 4$. Let $\{e_{ij}\mid 1\le 
i<j\le 4\}$ be the standard basis vectors of $\R^6$. By abuse of notation
$e_{ji}$ will stand for $-e_{ij}$ if $1\le i<j \le 4$.

Let $\rho_{S_4}$, the action of $S_4$ on $SO(3)$ be given 
as right multiplication by 
the rotation group of $U_3$, the
rhombic dodecahedron, which is easily seen to be isomorphic to the rotation
group of the cube. Moreover define $\tau_6$ to be the $S_4$ action on $\R^6$ 
given on the standard basis of $\R^6$ by the following rule: if 
$\sigma\in S_4$ is a permutation of the letters $1,2,3,4$ then 
$$\tau _6 (\sigma)(e_{ij})=sign(\sigma)e_{\sigma(i)\sigma(j)}.$$ 

The construction of the map $\Phi$ and the oddness of $F$ imply that 
$\Phi$ has to be an $S_4$ equivariant map from $(SO(3),\rho)$ to
$(\R^6,\tau_6)$. 

As we have seen above, to prove the theorem we need $A\in SO(3)$
such that $\Phi (A)$ lies in the $3$-subspace 
\begin{eqnarray*} V&:=&\{( \langle x, A(v_j-v_i) 
\rangle )_{1 \le i < j \le 4} \mid x \in {\R}^3 \}\\ &=&\{ (\langle A^*x,
v_j-v_i \rangle )_{1 \le i < j \le 4} \mid x \in {\R}^3 \}=\{ 
(\langle
y, v_j-v_i
\rangle )_{1 \le i < j \le 4} \mid y \in {\R}^3 \}
\end{eqnarray*}
of ${\R}^{6}$.

We prove more, namely that for any $S_4$ equivariant map 
$\Phi:(SO(3),\rho)\rightarrow (\R^6,\tau_6)$ there exists an 
$A\in SO(3)$ such that $\Phi(A)\in V$. 

The $3$-space $V$ is spanned by any three of the four vectors 
$e_{12}+e_{13}+e_{14},e_{21}+e_{23}+e_{24},e_{31}+e_{32}+e_{34},
e_{41}+e_{42}+e_{43}$, as we get these when $y=-v_1,-v_2,-v_3,-v_4$ 
respectively. Note that $V$ is an invariant subspace of $(\R^6,\tau_6)$. 
Let the 3-subspace $W$ of $\R^6$ be spanned by any three of the four vectors 
$e_{23}+e_{34}+e_{42},e_{31}+e_{14}+e_{43},e_{12}+e_{24}+e_{41},
e_{21}+e_{13}+e_{32}.$ Now one checks that 
$W$ is the orthogonal complement of
$V$ with respect to the standard Euclidean scalar product of $\R^6$. As the
action $\tau_6$ preserves this standard scalar product, it follows 
that $W$, being
the orthogonal complement of an invariant subspace $V$, is itself invariant 
under the action $\tau_6$. Let $\tau_W$ denote the $S_4$ action 
$\tau_6|_W$ on $W$. Now finding $A\in SO(3)$ with $\Phi(A)\in V$ is
equivalent to showing that the $S_4$ equivariant map 
$pr_W \Phi: (SO(3),\rho)\rightarrow (W,\tau_W)$ vanishes somewhere. 
By  Proposition~\ref{s4}
this is the case as $(W,\tau_W)$ is isomorphic to $(\R^3,\tau)$, which
is the action where $S_4$ acts as the symmetry group of the regular 
tetrahedron. This last statement can be seen by checking that $\tau_W$ 
faithfully permutes the four vectors 
$e_{23}+e_{34}+e_{42},e_{31}+e_{14}+e_{43},
e_{12}+e_{24}+e_{41},e_{21}+e_{13}+e_{32}$, which form a regular tetrahedron
in $W$. 

The result follows. 
\end{proof}

\par{\it Proof of Theorem~\ref{rhombic}.} As we have explained above,
 Theorem~\ref{concurrent} implies Makeev's conjecture for $n=3$, 
i.e., that $U_3$ is a universal cover in $\R^3$. Moreover $U_3$ is
the intersection of the $6$ strips corresponding to the $6$ edges of a
regular tetrahedron of edge length $1$, which is a rhombic
dodecahedron with distance of opposite faces equal to $1$.
$\square$  

\vskip.4cm 

A frequent application of universal covers is in the 
so-called Borsuk problem \cite{borsuk}: if $X \subset {\R}^n$ 
has diameter 
$1$, can it be decomposed into $n+1$ sets $X_1,...,X_{n+1}$ of smaller
diameters? 

We note that for all sufficiently large $n$ Borsuk's 
problem has a negative solution \cite{kahn-kalai} -- even for finite sets X --,
but the smallest $n$, for which a counterexample is known, is $n=561$, cf. 
\cite{boltyanski-martini-soltan} pp.209-226, and  
\cite{aigner-ziegler}, pp. 83-88. 

However for low dimensions the problem has an affirmative solution.

For $n=2$ the sharp answer is that one can even guarantee 
${\rm diam}(X_i) \le \sqrt 3 /2$ (\cite{gale}). 
The proof goes by applying P\'al's
theorem to cover $X$ with a regular hexagon $U_2$
with distance of opposite
sides $1$ and then cut $U_2$ into three congruent 
pentagons with diameter $\sqrt 3/2$.

For $n=3$ the positive answer has been 
proved first by \cite{eggleston2} and then 
\cite{heppes} and  \cite{grunbaum}. 
Heppes and Gr\"unbaum used for their proof a universal cover in $\R^3$, 
the regular octahedron
$O_3$, with distance of opposite faces equal to $1$, then chopped off
three 
vertices, still obtaining a universal cover. Then they decomposed this 
last set to four parts of diameters $0.998$ and $0.989$, respectively. 

For a nice exposition of the story of Borsuk's problem cf. 
\cite{boltyanski-martini-soltan}, pp. 209-226.

The universal cover $U_3$, given by Theorem~\ref{rhombic}, is 
intuitively smaller than $O_3$ (e.g., ${\rm diam}(O_3) =\sqrt 3,\,\,
{\rm diam} (U_3) =\sqrt 2$), but on the other hand is combinatorially more 
complex; so it is conceivable that with more amount of work one could 
substantially reduce the upper bound for ${\rm diam}(X_i)$. As we have been informed recently, \cite{makeev4} 
 already observed this point, moreover
referred to work of A. Evdokimov, who in this way sharpened the results of
Heppes and Gr\"unbaum, to 0.98.



\vskip.4cm

\noindent  Department of Mathematics, University of California
at Berkeley, Berkeley CA 94720-3840, USA; hausel@math.berkeley.edu 
\newline Mathematical Institute of the Hungarian Academy of Sciences,
H-1364 Budapest, Pf. 127, HUNGARY; makai@renyi.hu,
\newline L. E\"otv\"os University, Department of Analysis, H-1088 
M\'uzeum krt. 6-8, Budapest, HUNGARY; szucsandras@ludens.elte.hu

\end{document}